\ifdefined\pdfoutput\pdfoutput=1\fi % arXiv: pdfLaTeX output; ignored by XeLaTeX/LuaLaTeX
\documentclass[11pt]{article}
\usepackage[margin=1in]{geometry}
\usepackage[T1]{fontenc}
\usepackage{lmodern}
\usepackage{amsmath,amssymb,amsthm}
\usepackage{authblk}
\usepackage{xcolor}
\usepackage{blkarray}
\usepackage{dsfont}
\usepackage{bm}
\usepackage{amsmath}
\usepackage{amssymb}
\usepackage{graphicx}
\usepackage{subcaption}
\usepackage{tikz}
\usetikzlibrary{decorations.pathmorphing,patterns}
\usetikzlibrary{fit,shapes,arrows}

\newcommand{\tlam}{\tilde{\bar{\bm \lambda}}}
\newcommand{\tmu}{\tilde{\bar{\bm \mu}}}
\newcommand{\bx}{\bm{\bar{x}}}
\newcommand{\bu}{\bm{\bar{u}}}
\newcommand{\blam}{\bm{\bar{\lambda}}}
\newcommand{\bmu}{\bm{\bar{\mu}}}
\newcommand{\IR}{\mathbb{R}}

\newtheorem{thm}{\textbf{Theorem}}
\newtheorem{prob}{\textbf{Problem}}
\newtheorem{prop}{\textbf{Property}}
\newtheorem{lem}{\textbf{Lemma}}

\newtheorem{cor}{\textbf{Corollary}}
\newtheorem{asm}{\textbf{Assumption}}
\newtheorem{rem}{\textbf{Remark}}

\DeclareSymbolFont{extraup}{U}{zavm}{m}{n}
\DeclareMathSymbol{\vardiamond}{\mathalpha}{extraup}{87}
\usepackage{cite} % sorts and compresses numeric citations, e.g. [56-60]
\usepackage[colorlinks=true,linkcolor=blue,citecolor=blue,urlcolor=blue]{hyperref}

\title{Costate Estimation and Covector Mapping for Bernstein Polynomial Direct Methods in Optimal Control}

\author[1]{Venanzio Cichella\thanks{Corresponding author: venanzio-cichella@uiowa.edu}}
\author[2]{Isaac Kaminer}
\author[3]{Claire Walton}
\author[4]{Naira Hovakimyan}
\author[5]{Ant{\'{o}}nio M. Pascoal}
\affil[1]{Department of Mechanical Engineering, University of Iowa, Iowa City, IA 52242, USA}
\affil[2]{Department of Mechanical and Aerospace Engineering, Naval Postgraduate School, Monterey, CA 93940, USA}
\affil[3]{Department of Electrical and Computer Engineering, University of Texas at San Antonio, San Antonio, TX 78249, USA}
\affil[4]{Department of Mechanical Science and Engineering, University of Illinois at Urbana-Champaign, Urbana, IL 61801, USA}
\affil[5]{Institute for Systems and Robotics (ISR), Instituto Superior T\'{e}cnico (IST), University of Lisbon, 1049-001 Lisbon, Portugal}
\date{}

\begin{document}
\maketitle

\begin{abstract}

Bernstein polynomial approximation of continuous functions has a slower rate of convergence as compared to other approximation methods. \emph{``The fact seems to have precluded any numerical application of Bernstein polynomials from having been made. Perhaps they
will find application when the properties of the approximant in the large are of more importance than the closeness of the approximation.''} -- has remarked P.J. Davis in his 1963 book Interpolation and Approximation.
This paper presents a direct approximation method for nonlinear optimal control problems with mixed input and state constraints based on Bernstein polynomial approximation, in which the dynamic and path constraints are imposed at the Bernstein coefficients (control points). We provide a rigorous analysis showing that, under convergence assumptions on the control polygons, the proposed method yields consistent approximations of time continuous optimal control problems and can be used for costate estimation of the optimal control problems. In particular, multipliers satisfying a relaxed Karush–Kuhn–Tucker system of the resulting nonlinear program, scaled by the quadrature weight, approximate the costate and the multipliers of the optimal control problem. Owing to the structure of the Bernstein differentiation matrix, the discrete transversality conditions are contained in the relaxed Karush–Kuhn–Tucker system, and no additional closure conditions need to be added to it. This result leads to the formulation of the Covector Mapping Theorem for Bernstein polynomial approximation.
\end{abstract}


\section{Introduction}
\label{sec:Introduction}
{Motion planning plays an important role in enabling robotic systems accomplish tasks assigned to them autonomously, safely and reliably. Over the past decades many approaches to generating trajectories have been proposed. Examples include bug algorithms, artificial potential functions, roadmap path planners, cell decomposition
methods, and optimal control based trajectory generation. The reader is referred to \cite{ng2007performance,khatib1986real,siegwart2011introduction,choset2005principles,latombe2012robot,betts1998survey,von1992direct,lavalle2006planning} and references therein for detailed discussions and comparisons of these methods. Each technique has different advantages and disadvantages, and is best-suited for certain types of problems. Motion planning based on optimal control –i.e. optimal motion
planning– is particularly suitable for applications that require the trajectory to optimize some cost while guaranteeing satisfaction of a complex set of vehicle and problem constraints. These applications include multi-robot road search \cite{cichella2018cooperative}, coordinated tracking \cite{cichella2014faultynetworks}, optimal and constrained formation control \cite{sun2016optimal}, and adversarial swarm defense \cite{walton2022defense}.}

Optimal control problems that arise from robotics and motion planning applications are in general very complex. Finding a closed-form solution to these problems can be difficult or even impossible, and therefore they must be solved numerically. 
Numerical methods include indirect and direct methods \cite{rao2009survey}. Indirect methods solve the problems by converting them into boundary value problems. Then, the solutions are found by solving systems of differential equations. On the other hand, direct methods are based on transcribing optimal control problems into nonlinear programming problems (NLPs) using some discretization scheme \cite{rao2009survey,betts2010practical,betts1998survey,conway2012survey}. These NLPs can be solved using ready-to-use NLP solvers (e.g. MATLAB, SNOPT, etc.) and do not require calculation of costate and adjoint variables as indirect methods do. 

A wide range of direct methods that uses different discretization schemes have been developed, including direct single shootings, direct multiple shooting, and direct collocation methods \cite{rao2009survey,betts2010practical,betts1998survey,conway2012survey}.  
The software packages that implement some of these methods (e.g., PSOPT \cite{becerra2010solving}, NLOptControl \cite{febbo2020nloptcontrol}, GPOPS II \cite{patterson2014gpops}, PROPT \cite{rutquist2010propt}, DIDO \cite{ross2020enhancements}, and CasADi \cite{andersson2019casadi}) are particularly relevant; some of these have been
applied successfully to solve a wide range of real-world problems \cite{fahroo2006discrete,bollino2007pseudospectral,gong2009pseudospectral,bedrossian2009zero,bollino2008collision,bedrossian2007first,ross2012review}.
{Theoretical results in the literature on direct methods include that on consistency of approximation theory, see \cite{elijah1997optimization}, which provides a framework to assess the convergence properties of Euler and Runge-Kutta discretization schemes. Motivated by the consistency of approximation theory, direct methods that use different discretization schemes have been developed, including Pseudospectral methods based on Legendre, Chebyshev, and Lagrange polynomials \cite{ross2012review}. }
% Relevant works in the literature on direct methods include that of of Polak on consistency of approximation theory, see \cite{elijah1997optimization}, which provides a theoretical framework to assess the convergence properties of direct methods. Motivated by the consistency of approximation theory, direct methods that use different discretization schemes have been developed, including Euler \cite{elijah1997optimization}, Runge-Kutta \cite{schwartz1996consistent}, and Pseudospectral \cite{ross2012review} methods. 
One drawback of direct methods is that the costate of the original optimal control problem can not be readily obtained from the approximated solution. {Nevertheless, in several applications -- such as motion planning and control for safety-critical robotic systems -- the knowledge of the costate is important because it allows for the evaluation of the fulfillment of necessary conditions of optimality. This evaluation, in turn, provides important insights on the validity and optimality of the solution.} Therefore, approaches for obtaining estimates of the costate from direct methods have been proposed in the literature on direct collocation \cite{fahroo2001costate,darby2011costate,hager2000runge} and direct shooting \cite{grimm1997adjoint}.

Over the past several years, our group has developed a family of methods built on Bernstein polynomials, motivated by the observation that their geometric properties make them well suited to describe and constrain trajectories. The first applications were offline optimal motion planning problems. In \cite{cichella2018optimal,cichella2020optimal} we showed that the end-point, convex hull and de Casteljau properties of Bernstein polynomials allow for efficient and conservative enforcement of state and input constraints, such as collision avoidance between vehicles, and these ideas were later extended to composite (piecewise) Bernstein approximants \cite{maclin2024optimal} and to trajectory generation for autonomous aircraft operating near airports \cite{maclin2025transformer}. The resulting transcription methods were collected in the software toolkit BeBOT \cite{kielas2019bebot,kielas2022bernstein}. Because the same properties lead to small and well-structured nonlinear programs, Bernstein-based formulations have also been used online: for persistent and randomized monitoring \cite{stephens2019randomized,kielas2021persistent}, for motion planning aimed at target localization \cite{tabasso2021optimal,tabasso2024online}, in combination with learned warm starts \cite{maclin2025realtime}, for geodesic trajectory planning \cite{gorman2026solving}, for control-agnostic model predictive planning \cite{gorman2027model}, and for model predictive control of autonomous underwater vehicles \cite{hammond2025adaptive,petrov2026adaptive}. Their bivariate counterparts, Bernstein surfaces, have been used to plan the motion of soft continuum robots and of systems governed by partial differential equations \cite{hammond2024path,hammond2025constrained,golestaneh2024scalable,hammond2025scalable}. More recently, Bernstein (B\'{e}zier) signals have been used as models of exogenous signals in output regulation based on the internal model principle \cite{cichella2026robust,mimmo2026tracking,cichella2026polynomial}. Additional works that exploit the properties of Bernstein polynomials for nonlinear optimal control can be found in \cite{ricciardi2018direct,choe2016cooperative,ghomanjani2012optimal,huo2022fast,zhao2017split}.

Across these applications, a recurring question is whether the solutions of the discretized problems are faithful approximations of the solutions of the underlying continuous-time problems. This question was addressed in \cite{cichella2018bernstein,cichellaBernsteinCDC,cichella2020optimal}, where we established feasibility and consistency results for the primal problem, and in \cite{cichella2022consistency}, where we extended the analysis to the dual problem and to costate estimation for a discretization in which the dynamic and path constraints are evaluated at time nodes along the Bernstein polynomials. The present paper revisits the latter analysis for the formulation of \cite{cichella2020optimal}, in which the dynamic and path constraints are imposed at the Bernstein coefficients (control points). This requires new proofs, since the discrete multipliers are now associated with the control points rather than with time nodes, and it leads to a simpler set of discrete optimality conditions.
 In particular, we propose an approach to approximate the costate of a general non-linear optimal control problem of Bolza type using multipliers satisfying a relaxed Karush–Kuhn–Tucker system of the Bernstein polynomial-based discrete approximation. We derive transformations that relate the Lagrange multipliers of the non-linear programming problem to the costate of the original optimal control problem. These transformations are often referred to as \emph{covector mapping} {in the literature on direct methods for optimal control \cite{ross2012review,elijah1997optimization,schwartz1996consistent}. Finally, we demonstrate uniform convergence properties of the method.}    
 
The paper is structured as follows: in Section~\ref{sec:MathematicalPreliminaries} we present the notation and the mathematical results, which will be used later in the paper. Section~\ref{sec:ProblemFormulation} introduces the optimal control problem of interest and some related assumptions, and presents the approximation method based on Bernstein approximation that approximates the optimal control problem into an NLP. In Section \ref{sec:costate} we derive the relaxed KKT conditions associated with the NLP. Section \ref{sec:costateresults} compares these conditions to the first order optimality conditions for the original optimal control problem and states the Covector Mapping Theorem for Bernstein approximation. Numerical examples are discussed in Section~\ref{sec:numericalresults},  The paper ends with conclusions in Section~\ref{sec:conclusions}.
%%%%%%%%%%%%%%%%%%%%%%%%%%%%%%%%%%%%%%%%%

%============== NOTATION =========================
%=================================================

\section{Notation and mathematical background}\label{sec:MathematicalPreliminaries}
Vector valued functions are denoted by bold letters, $\bm{x}(t) =[x_1(t) \, , \, \ldots \, , \, x_n(t)]^\top$, while vectors are denoted by bold letters with an upper bar, $\bx = [x_1 \, , \, \ldots \, , \, x_n]^\top \in \IR^{n}$. The symbol $\mathcal{C}^r$ denotes the space of functions with $r$ continuous derivatives. $\mathcal{C}^r_n$ denotes the space of $n$-vector valued functions in $\mathcal{C}^r$. $||\cdot||$ denotes the Euclidean norm, $||\bx|| = \sqrt{x_1^2 + \ldots + x_n^2}$.

The Bernstein basis polynomials of degree $N$ are defined as
$$ b_{j,N}(t) = \binom{N}{j} t^j(1-t)^{N-j}  \, , \qquad t \in [0,1] \, ,$$
for $j=0,\ldots,N$, with
$
\binom{N}{j} = \frac{N !}{j!(N-j)!}\,.
$
%They were originally introduced by the mathematician Sergei Natanovich Bernstein in 1912 to facilitate a constructive proof of the Weierstrass approximation theorem \cite{bernstein1912}. 
An $N$th order Bernstein polynomial $x_N:[0,1]\to\IR$ is a linear combination of $N+1$ Bernstein basis polynomials of order $N$, i.e.
$$
x_N(t) = \sum_{j=0}^N \bar{x}_j b_{j,N}(t) \, , \qquad t \in [0,1] \, ,
$$
where $\bar{x}_j \in \IR$, $j=0,\ldots,N$, are referred to as Bernstein coefficients. For the sake of generality, and with a slight abuse of terminology, in this paper we extend the definition of a Bernstein polynomial given above to a vector of $N$th order polynomials $\bm{x}_N : [0,1] \to \IR^n$ expressed in the following form
\begin{equation} \label{eq:beziercurve}
\bm{x}_N(t) = \sum_{j=0}^N  \bx_{j,N} b_{j,N}(t)  \, , \qquad t \in [0,1] \, ,
\end{equation}
where $\bx_{0,N},\ldots,\bx_{N,N} \in \mathbb{R}^n$. 
The Bernstein coefficients are also referred to as \emph{control points}. Throughout the paper, $t_j = \frac{j}{N}$, $j=0,\ldots,N$, denote equidistant time nodes. The \emph{control polygon} $\bm{X}_N:[0,1]\to\IR^n$ of $\bm{x}_N(t)$ is the piecewise linear function obtained by connecting the Bernstein coefficients in order, i.e.,
\begin{equation} \label{eq:controlpolygon}
\bm{X}_N(t) := \bx_{j,N}\frac{t_{j+1}-t}{t_{j+1}-t_j} + \bx_{j+1,N}\frac{t-t_j}{t_{j+1}-t_j} \, , \qquad t\in[t_j,t_{j+1}] \, , \quad j=0,\ldots,N-1 \, ,
\end{equation}
and its slopes are the difference quotients
\begin{equation} \label{eq:slopes}
\bm{\Delta}_{i}(\bx_N) := N\left(\bx_{i+1,N}-\bx_{i,N}\right) \, , \qquad i = 0,\ldots,N-1 \, .
\end{equation}

%Bernstein polynomials were popularized by Pierre B\'{e}zier in the early 1960s as useful tools for geometric design (B\'{e}zier used Bernstein polynomials to design the shape of the cars at the Renault company in France), and are now widely used in computer graphics, animations and type fonts such as postscript fonts and true type fonts. For this reason, the Bernstein polynomial introduced in Equation \eqref{eq:beziercurve} is often referred to as a B{\' e}zier curve, especially when used to describe a spatial curve.

%Bernstein polynomials possess favorable geometric and numerical properties that can be exploited in many application domains.
In what follows, we provide a review of numerical properties of Bernstein polynomials that are used throughout this paper.
%For an extensive review on Bernstein polynomials and their properties the reader is referred to \cite{farouki2012bernstein}.
The derivative and integral of a Bernstein polynomial $\bm{x}_N(t)$ can be easily computed as
\begin{equation*}
\dot{\bm{x}}_N(t) = N \sum_{j = 0}^{N-1} (\bx_{j+1,N}-\bx_{j,N}) b_{j,N-1}(t) \,
\end{equation*}
and
\begin{equation} \label{eq:integralbezier}
\int_0^1{\bm{x}}_N(t)dt = w \sum_{j = 0}^{N} \bx_{j,N} \, , \qquad w = \frac{1}{N+1} \, ,
\end{equation}
respectively.
In other words, the slopes $\bm{\Delta}_{j}(\bx_N)$ are the Bernstein coefficients of $\dot{\bm{x}}_N(t)$ in the basis of degree $N-1$. Elevating their degree from $N-1$ to $N$ \cite{farouki2012bernstein}, the derivative can be expressed exactly in the basis of degree $N$:
\begin{equation} \label{eq:derivativeDN}
\dot{\bm{x}}_N(t) = \sum_{j=0}^N \Big( \sum_{k=0}^N \bx_{k,N} D_{k,j} \Big) b_{j,N}(t) \, , \qquad t \in [0,1] \, ,
\end{equation}
where $D_{k,j}$ is the $(k,j)$th entry of the square differentiation matrix $\bm{D}_N \in \IR^{(N+1)\times(N+1)}$, which is tridiagonal with entries
\begin{equation} \label{eq:DNentries}
D_{k,k} = 2k-N \, , \qquad D_{k,k-1} = N-k+1 \, , \qquad D_{k,k+1} = -(k+1) \, ,
\end{equation}
$k=0,\ldots,N$, and zero otherwise.
Bernstein polynomials can be used to approximate smooth functions.
Consider a $n$-vector valued function $\bm{x}:[0,1] \to \IR^n$. The $N$th order \emph{Bernstein approximation} of $\bm{x}(t)$ is a vector of Bernstein polynomials $\bm{x}_N (t)$ computed as in \eqref{eq:beziercurve} with $\bx_{j,N} = \bm{x}(t_j)$ and $t_j = \frac{j}{N}$ for all $j = 0,\ldots,N$. Namely,
\begin{equation} \label{eq:bernapprox}
{ \bm{x}(t) \approx } \bm{x}_N(t) = \sum_{j=0}^N  \bm{x}(t_j) b_{j,N}(t)  \, , \qquad t_j = \frac{j}{N} \, .
\end{equation}
The following results hold for Bernstein approximations.
\begin{lem}[Uniform convergence of Bernstein approximation]
	\label{lem:bernsteinapprox}
	Let $\bm{x}(t) \in \mathcal{C}^{0}_n$ on $[0,1]$, and let $\bm{x}_N(t)$ be computed as in Equation \eqref{eq:bernapprox}. Then, for arbitrary order of approximation $N \in \mathbb{Z}^+$, the Bernstein approximation $\bm{x}_N(t)$ satisfies
	\begin{equation*}
	||\bm{x}_N(t)-\bm{x}(t)|| \leq C_0 W_x(N^{-\frac{1}{2}}) \, ,
	\end{equation*}
	where $C_0$ is a positive constant satisfying $C_0 < 5n/4$, and $W_x(\cdot)$ is the modulus of continuity of $\bm{x}(t)$ in $[0,1]$ \cite{bojanic1989rate,popoviciu1935approximation,sikkema1961wert}. 
	%In other words, the Bernstein approximation $\bm{x}_N(t)$ converges uniformly to $\bm{x}(t)$ on $[0,1]$. 
	Moreover, if $\bm{x}(t) \in \mathcal{C}^{1}_n$, then
	\begin{equation*}
	\Vert \dot{\bm{x}}_N(t)-\dot{\bm{x}}(t) \Vert \leq C_1 W_{x^\prime}(N^{-\frac{1}{2}}) \, ,
	\end{equation*}
	where $C_1$ is a positive constant satisfying $C_1 < 9n/4$ and $W_{x^{\prime}}(\cdot)$ is the modulus of continuity of $\dot{\bm{x}}(t)$ in $[0,1]$ \cite{powell1981approximation}.
\end{lem}
\begin{lem} \cite{floater2005convergence} \label{lem:bernsteinapproxC2}
Assume $\bm{x}(t) \in \mathcal{C}^{r+2}_n$, $r\geq 0$, and let $\bm{x}_N(t)$ be computed as in Equation \eqref{eq:bernapprox}. Let $\bm{x}^{(r)}(t)$ denote the $r$th derivative of $\bm{x}(t)$. Then, the following inequalities hold for all $t\in [0,1]$:
\begin{equation*}
\begin{split}
||\bm{x}_N(t)-\bm{x}(t) || & \leq \frac{C_0}{N} \, ,\\
%||\bm{f}^\prime_N(t)-\bm{f}^\prime(t) || & \leq \frac{C_1}{N} \\
& \, \, \, \vdots \\
||\bm{x}^{(r)}_N(t)-\bm{x}^{(r)}(t) || & \leq \frac{C_r}{N} \, ,
\end{split}
\end{equation*}
 where $C_0,\ldots,C_r$ are independent of the order of approximation $N$.
\end{lem}

\begin{lem} \label{lem:quadrature}
If $\bm{x}(t) \in \mathcal{C}^{0}_n$ on $[0,1]$, then we have
$$
\left\Vert \int_0^1 \bm{x}(t) dt - w \sum_{j=0}^N \bm{x}\left(  \frac{j}{N} \right) \right\Vert \leq C_I W_x(N^{-\frac{1}{2}}) \, ,
$$
with $w = \frac{1}{N+1}$, where $C_I> 0 $ is independent of $N$.
Moreover, if $\bm{x}(t) \in \mathcal{C}^{2}_n$, then
$$
\left\Vert \int_0^1 \bm{x}(t) dt - w \sum_{j=0}^N \bm{x}\left(  \frac{j}{N} \right) \right\Vert \leq \frac{C_I}{N} \, .
$$
\end{lem}
The Lemma above follows directly from Lemmas \ref{lem:bernsteinapprox} and \ref{lem:bernsteinapproxC2} and Equation \eqref{eq:integralbezier}.
In the above lemmas, the modulus of continuity of $\bm{x}:[0,1]\to\IR^n$ is
\begin{equation} \label{eq:modulus}
W_x(\delta) = \max_{|t_1-t_2|\leq \delta} ||\bm{x}(t_1)-\bm{x}(t_2)|| \, , \qquad t_1,t_2\in[0,1] \, .
\end{equation}
If $\bm{x}\in\mathcal{C}^0_n$, then $\lim_{\delta\to 0}W_x(\delta)=0$, and $W_x(\delta) \leq \delta ||\dot{\bm{x}}||_\infty$ if $\bm{x}\in\mathcal{C}^1_n$, where, for a continuous function, $||\bm{x}||_\infty := \max_{t\in[0,1]}||\bm{x}(t)||$.

\begin{lem} \label{lem:controlpolygon}
Let $\bm{x}(t)\in\mathcal{C}^0_n$ on $[0,1]$, let $\bm{x}_N(t)$ be a Bernstein polynomial with coefficients $\bx_{j,N}$, and let $\bm{X}_N(t)$ be its control polygon. Then
\begin{align*}
\lim_{N\to\infty}\max_{t\in[0,1]}||\bm{X}_N(t)-\bm{x}(t)||=0 & \implies \lim_{N\to\infty}\max_{j}||\bx_{j,N}-\bm{x}(t_j)||=0 \\
& \implies \lim_{N\to\infty}\max_{t\in[0,1]}||\bm{x}_N(t)-\bm{x}(t)||=0 \, .
\end{align*}
\end{lem}
The first implication follows from $\bm{X}_N(t_j)=\bx_{j,N}$. For the second, $||\bm{x}_N(t)-\bm{x}(t)|| \leq \max_j||\bx_{j,N}-\bm{x}(t_j)|| + ||\sum_j \bm{x}(t_j)b_{j,N}(t)-\bm{x}(t)||$, since the basis polynomials are nonnegative and sum to one; the last term vanishes uniformly by Lemma \ref{lem:bernsteinapprox}.

The next lemma characterizes the action of the differentiation matrix $\bm{D}_N$ on a sequence of Bernstein coefficients, both from the right (as in Equation \eqref{eq:derivativeDN}) and from the left. The former is used to approximate the dynamics of Problem $P$, while the latter appears in the KKT conditions of its discretization and approximates the costate dynamics.
\begin{lem}[Differentiation matrix and its adjoint action] \label{lem:Dproperties}
Let $\bar{\bm{a}}_{0},\ldots,\bar{\bm{a}}_{N} \in\IR^n$, $N\geq 1$, let $\bm{\Delta}_i = N(\bar{\bm{a}}_{i+1}-\bar{\bm{a}}_{i})$, $i=0,\ldots,N-1$, let $w=\frac{1}{N+1}$, and let $D_{k,j}$ be given by Equation \eqref{eq:DNentries}. Then:
\begin{enumerate}
\item[(i)] for $j = 0,\ldots,N$,
$\displaystyle
\sum_{k=0}^N \bar{\bm{a}}_{k} D_{k,j} = \frac{j}{N}\bm{\Delta}_{j-1}+\frac{N-j}{N}\bm{\Delta}_{j} \, ,
$
where the terms $\bm{\Delta}_{-1}$ and $\bm{\Delta}_{N}$, which appear for $j=0$ and $j=N$, have zero weight;
\item[(ii)] for $k = 1,\ldots,N-1$,
$\displaystyle
-\sum_{j=0}^N D_{k,j} \bar{\bm{a}}_{j} = \frac{N-k+1}{N}\bm{\Delta}_{k-1}+\frac{k+1}{N}\bm{\Delta}_{k} \, ;
$
\item[(iii)]
$\displaystyle
-w\sum_{j=0}^N D_{0,j} \bar{\bm{a}}_{j} = \bar{\bm{a}}_{0} + \frac{w}{N}\bm{\Delta}_{0} \, , \qquad
-w\sum_{j=0}^N D_{N,j} \bar{\bm{a}}_{j} = -\bar{\bm{a}}_{N} + \frac{w}{N}\bm{\Delta}_{N-1} \, .
$
\end{enumerate}
Moreover, if $\bar{\bm{a}}_{j} = \bm{a}(t_j)$, $j=0,\ldots,N$, for some $\bm{a}(t)\in\mathcal{C}^1_n$ on $[0,1]$, then $||\bm{\Delta}_i|| \leq \sqrt{n}\,||\dot{\bm{a}}||_\infty$ and
\begin{align}
& \Big\Vert \sum_{k=0}^N \bar{\bm{a}}_{k} D_{k,j} - \dot{\bm{a}}(t_j) \Big\Vert \leq \sqrt{n} \, W_{a^\prime}\left(N^{-1}\right) \, , \qquad j=0,\ldots,N \, , \label{eq:Dright} \\
& \Big\Vert -\sum_{j=0}^N D_{k,j} \bar{\bm{a}}_{j} - \dot{\bm{a}}(t_k) \Big\Vert \leq \sqrt{n} \, W_{a^\prime}\left(N^{-1}\right) + \frac{2\sqrt{n}}{N}||\dot{\bm{a}}||_\infty \, , \qquad k=1,\ldots,N-1 \, , \label{eq:Dleft}
\end{align}
where $W_{a^\prime}(\cdot)$ is the modulus of continuity of $\dot{\bm{a}}(t)$.
\end{lem}
\noindent\textbf{Proof:} (i) By Equation \eqref{eq:DNentries}, $\sum_{k}\bar{\bm{a}}_{k} D_{k,j} = -j\bar{\bm{a}}_{j-1}+(2j-N)\bar{\bm{a}}_{j}+(N-j)\bar{\bm{a}}_{j+1} = j(\bar{\bm{a}}_{j}-\bar{\bm{a}}_{j-1})+(N-j)(\bar{\bm{a}}_{j+1}-\bar{\bm{a}}_{j})$, which gives (i).
(ii) For $1\leq k\leq N-1$, $-\sum_{j}D_{k,j}\bar{\bm{a}}_{j} = -(N-k+1)\bar{\bm{a}}_{k-1}-(2k-N)\bar{\bm{a}}_{k}+(k+1)\bar{\bm{a}}_{k+1} = (N-k+1)(\bar{\bm{a}}_{k}-\bar{\bm{a}}_{k-1})+(k+1)(\bar{\bm{a}}_{k+1}-\bar{\bm{a}}_{k})$, where we used $(N-k+1)+(2k-N)-(k+1)=0$.
(iii) For $k=0$, $-\sum_{j}D_{0,j}\bar{\bm{a}}_{j} = N\bar{\bm{a}}_{0}+\bar{\bm{a}}_{1} = (N+1)\bar{\bm{a}}_{0}+(\bar{\bm{a}}_{1}-\bar{\bm{a}}_{0})$; multiplying by $w=\frac{1}{N+1}$ yields the first identity. The second follows from $-\sum_{j}D_{N,j}\bar{\bm{a}}_{j} = -\bar{\bm{a}}_{N-1}-N\bar{\bm{a}}_{N} = -(N+1)\bar{\bm{a}}_{N}+(\bar{\bm{a}}_{N}-\bar{\bm{a}}_{N-1})$.
Now let $\bar{\bm{a}}_{j}=\bm{a}(t_j)$. By the mean value theorem applied to each component $a_l(t)$, $l=1,\ldots,n$, the $l$th component of $\bm{\Delta}_i$ equals $\dot{a}_l(\zeta_{i,l})$ for some $\zeta_{i,l}\in[t_i,t_{i+1}]$, which gives $||\bm{\Delta}_i||\leq\sqrt{n}\,||\dot{\bm{a}}||_\infty$. In (i), the weights are nonnegative, sum to one, and multiply slopes on the two segments adjacent to $t_j$, so that $|\zeta_{i,l}-t_j|\leq N^{-1}$; each component of the left hand side of \eqref{eq:Dright} is therefore bounded by $W_{a^\prime}(N^{-1})$, which gives \eqref{eq:Dright}. In (ii), let $\alpha_k = \frac{N-k+1}{N+2}\in(0,1)$; then $-\sum_{j}D_{k,j}\bar{\bm{a}}_{j} = \frac{N+2}{N}\big(\alpha_k\bm{\Delta}_{k-1}+(1-\alpha_k)\bm{\Delta}_{k}\big)$. By the same argument, each component of $\alpha_k\bm{\Delta}_{k-1}+(1-\alpha_k)\bm{\Delta}_{k}-\dot{\bm{a}}(t_k)$ is bounded by $W_{a^\prime}(N^{-1})$, while $||\frac{2}{N}(\alpha_k\bm{\Delta}_{k-1}+(1-\alpha_k)\bm{\Delta}_{k})||\leq\frac{2\sqrt{n}}{N}||\dot{\bm{a}}||_\infty$. The triangle inequality yields \eqref{eq:Dleft}.

\begin{rem} \label{rem:adjoint}
Lemma \ref{lem:Dproperties} shows that, acting from the right, $\bm{D}_N$ maps the samples of a $\mathcal{C}^1$ function to a convex combination of the adjacent slopes of its control polygon, i.e., to an approximation of its derivative. Acting from the left, the interior rows of $\bm{D}_N$ map the samples to $-\frac{N+2}{N}$ times a convex combination of the adjacent slopes, i.e., to an approximation of \emph{minus} the derivative, whereas the first and last rows contain the terms $\frac{1}{w}\bar{\bm{a}}_0$ and $-\frac{1}{w}\bar{\bm{a}}_N$. These boundary terms are the discrete counterpart of the boundary terms of an integration by parts, since $\sum_{j}\bar{\bm{b}}_j^\top\sum_k\bar{\bm{a}}_kD_{k,j} = \sum_k\bar{\bm{a}}_k^\top\sum_jD_{k,j}\bar{\bm{b}}_j$ for any $\bar{\bm{b}}_0,\ldots,\bar{\bm{b}}_N\in\IR^n$. As shown in Section \ref{sec:costate}, they produce the transversality conditions of the costate.
\end{rem}

The following property of Bernstein polynomials is relevant to this paper.
\begin{prop}[End point values] \label{app.prop:endpoint}
The Bernstein polynomial given by Equation \eqref{eq:beziercurve} satisfies $\bm{x}_N (0) = \bm{\bar{x}}_{0,N}$ and $\bm{x}_N(1) = \bm{\bar{x}}_{N,N}$.
\end{prop}
\begin{prop}[Convex hull] \label{app.prop:convexhull}
A Bernstein polynomial is completely contained in the convex hull of its Bernstein coefficients. In particular, for a scalar Bernstein polynomial, $\min_j \bar{x}_{j,N} \leq x_N(t) \leq \max_j \bar{x}_{j,N}$ for all $t\in[0,1]$.
\end{prop}

\section{Problem formulation} \label{sec:ProblemFormulation}
This paper considers the following optimal control problem:
\begin{prob}[Problem $P$]
Determine $\bm{x}: [0,1] \to \IR^{n_x}$ and $\bm{u}: [0,1] \to \IR^{n_u}$ that minimize
\begin{equation} \label{eq:costfunc}
I(\bm{x}(t),\bm{u}(t)) = E(\bm{x}(0),\bm{x}(1))+ \int_0^{1} F(\bm{x}(t),\bm{u}(t))dt \,,
\end{equation}
subject to
\begin{align}
& \dot{\bm{x}} = \bm{f}(\bm{x}(t),\bm{u}(t))\, , \quad \forall t \in[0,1], \label{eq:dynamicconstraint} \\
& \bm{e}(\bm{x}(0),\bm{x}(1)) = \bm{0} \, , \label{eq:equalityconstraint} \\
& \bm{h}(\bm{x}(t),\bm{u}(t)) \leq \bm{0} \, , \quad \forall t\in [0,1] \label{eq:inequalityconstraint} \, ,
\end{align}
where $E:\IR^{n_x}\times \IR^{n_x}\to \IR$ and $F:\IR^{n_x}\times \IR^{n_u}\to \IR$ are the terminal and running costs, respectively, $\bm{f}:\IR^{n_x}\times \IR^{n_u}\to \IR^{n_x}$ describes the system dynamics, $\bm{e}:\IR^{n_x}\times \IR^{n_x}\to \IR^{n_e}$ is the vector of boundary conditions, and $\bm{h}:\IR^{n_x} \times \IR^{n_u}\to \IR^{n_h}$ is the vector of state and input constraints.
\end{prob}

%The following assumptions are made:
%\begin{asm} \label{asm:Flip}
%$E$, $F$, $\bm{f}$, $\bm{e}$, and $\bm{h}$ are Lipschitz continuous with respect to their arguments.
%%$F(\cdot)\in \mathcal{C}^2$.
%\end{asm}
%\begin{asm} \label{asm:pexists}
%Problem $P$ admits optimal solutions $\bm{x}^*(t)$ and $\bm{u}^*(t)$ that satisfy $\bm{x}^*(t) \in \mathcal{C}^1_{n_x}$ and $\bm{u}^*(t) \in \mathcal{C}^0_{n_u}$.
%\end{asm}

%==============================================================================

%============ Approximation ===============
%================================================

% %\section{Bernstein approximation of Problem P} \label{sec:BernApprox}

Next, we formulate a discretized version of Problem $P$, here referred to as Problem $P_N$, where $N$ denotes the \emph{order of approximation}.
This requires that we approximate the input and state functions, the cost function, the system dynamics, and the equality and inequality constraints in Problem $P$.
First, consider the following $N$th order vectors of Bernstein polynomials:
\begin{equation} \label{eq:polappr}
\bm{x}_N(t) =   \sum_{j=0}^N \bx_{j,N} b_{j,N}(t) \, , \quad  \bm{u}_N(t) =   \sum_{j=0}^N \bu_{j,N} b_{j,N}(t) ,
\end{equation}
with $\bm{x}_N: [0,1] \to \IR^{n_x}$, $\bm{u}_N: [0,1] \to \IR^{n_u}$, $\bx_{j,N} \in \IR^{n_x}$ and $\bu_{j,N} \in \IR^{n_u}$. Let $\bx_N \in \IR^{n_x \times (N+1)}$ and $\bu_N \in \IR^{n_u \times (N+1)}$ be defined as
$$\bx_N = [\bx_{0,N} \, , \, \ldots \, , \, \bx_{N,N}] , \qquad \bu_N = [\bu_{0,N} \, , \, \ldots \, , \, \bu_{N,N}].$$
Then Problem $P_N$ can be stated as follows:
\begin{prob}[Problem $P_N$]
Determine $\bx_N$ and $\bu_N$ that minimize
\begin{equation} \label{eq:discCVcostfunc}
I_N(\bx_N,\bu_N) =
E(\bx_{0,N},\bx_{N,N})+ w \sum_{j=0}^N F(\bx_{j,N},\bu_{j,N}) \, ,
\end{equation}
subject to
\begin{align}
& \Big\Vert \sum_{k=0}^N \bx_{k,N}D_{k,j} - \bm{f}(\bx_{j,N},\bu_{j,N})\Big\Vert \leq \delta_P^N \, ,
\qquad
\forall j=0,\ldots,N \,,
\label{eq:discdynamicconstraint} \\
& \bm{e}(\bx_{0,N},\bx_{N,N}) = \bm{0} \, , \label{eq:discequalityconstraint} \\
& \bm{h}(\bx_{j,N},\bu_{j,N}) \leq   \bm{0} \, , \qquad  \forall j=0,\ldots,N  \label{eq:discinequalityconstraint} \, ,
\end{align}
where $w=\frac{1}{N+1}$, $D_{k,j}$ is the $(k,j)$th entry of the differentiation matrix $\bm{D}_N$ in Equation \eqref{eq:DNentries}, and $\delta_P^N$ is a small positive number that depends on $N$ and satisfies $\lim_{N \to \infty} \delta_P^N = 0$.
\end{prob}
\vspace{0.25cm}

In Problem $P_N$, the integral term of the cost function is replaced by a weighted sum of the running cost evaluated at the Bernstein coefficients, and the end point cost and the boundary conditions are evaluated at the first and last coefficients, which coincide with the end points of $\bm{x}_N(t)$ (Property \ref{app.prop:endpoint}). The dynamic constraint and the state and input constraints are imposed at the Bernstein coefficients (control points), rather than at time nodes along the polynomials. By Equation \eqref{eq:derivativeDN}, $\sum_{k}\bx_{k,N}D_{k,j}$ is the $j$th Bernstein coefficient of $\dot{\bm{x}}_N(t)$ in the basis of degree $N$.
\vspace{0.25cm}

\begin{rem} \label{rem:relaxation}
Compared to the dynamic constraint of Problem $P$, the dynamic constraint given by Equation \eqref{eq:discdynamicconstraint} is \emph{relaxed}. Motivated by previous work on consistency of approximation theory \cite{elijah1997optimization}, the bound $\delta_P^N$, referred to as \emph{relaxation bound}, is introduced to guarantee that Problem $P_N$ has a feasible solution. As shown in Theorem \ref{thm:existence_dual} below, the relaxation bound can be made arbitrarily small by choosing a sufficiently large order of approximation $N$. The inequality constraint \eqref{eq:discinequalityconstraint} is not relaxed. If each component of $\bm{h}$ is convex, then, since $(\bm{x}_N(t),\bm{u}_N(t))$ lies in the convex hull of the pairs $(\bx_{j,N},\bu_{j,N})$ (Property \ref{app.prop:convexhull}), Equation \eqref{eq:discinequalityconstraint} implies $\bm{h}(\bm{x}_N(t),\bm{u}_N(t))\leq\bm{0}$ for all $t\in[0,1]$ and any $N$.
\end{rem}

%======================================================================
\begin{rem} 
The outcome of Problem $P_N$ is a set of optimal Bernstein coefficients $\bx^*_N$ and $\bu^*_N$ that determine the vectors of Bernstein polynomials $\bm{x}_N^*(t)$ and $\bm{u}_N^*(t)$, i.e.
\begin{equation} \label{eq:approxoptimal}
\bm{x}_N^* (t) = \sum_{j=0}^N \bx_{j,N}^* b_{j,N}(t) \, , \quad \bm{u}_N^* (t) = \sum_{j=0}^N \bu_{j,N}^* b_{j,N}(t)\, .
\end{equation}
In our previous work, see \cite{cichella2020optimal}, we provide theoretical results demonstrating:
(i) the existence of a feasible solution to Problem $P_{N}$,
and (ii) the convergence of the pair $(\bm{x}_N^* (t),\bm{u}_N^* (t))$ to the optimal solution of Problem $P$, given by $(\bm{x}^* (t),\bm{u}^* (t))$. Nevertheless, the present paper focuses on the existence and convergence of the estimates of the costates of Problem $P$, which are introduced next.
These results are established in \cite{cichella2020optimal} for the formulation adopted here, in which the dynamic and path constraints are imposed at the Bernstein coefficients. Primal feasibility and consistency are also recalled in Section \ref{sec:costateresults}, see Theorems \ref{thm:existence_dual} and \ref{thm:consistency_dual}.
\end{rem}

\section{Costate estimation for Problem $P$} \label{sec:costate}
\subsection{First order optimality conditions of Problem P}
We start by deriving the first order necessary conditions for Problem $P$. Let ${\bm \lambda}(t): [0,1] \to \IR^{n_x}$ be the costate trajectory, and let ${\bm \mu}(t): [0,1] \to \IR^{n_h}$ and ${\bm \nu}\in \IR^{n_e}$ be the multipliers.
By defining the Lagrangian of the Hamiltonian (also known as the D-form \cite{hartl1995survey}) as
\begin{equation*}
\begin{split}
& \mathcal{L}(\bm{x}(t),\bm{u}(t),\bm{\lambda}(t),\bm{\mu}(t)) = \\ & \qquad  \mathcal{H}(\bm{x}(t),\bm{u}(t),\bm{\lambda}(t)) + {\bm \mu}^\top (t) \bm{h}(\bm{x}(t),\bm{u}(t)) \, ,
\end{split}
\end{equation*}
where the Hamiltonian $\mathcal{H}$ is given by 
$$
\mathcal{H}(\bm{x}(t),\bm{u}(t),\bm{\lambda}(t)) = F({\bm x}(t),{\bm u}(t)) + {\bm \lambda}^\top(t) \bm{f}(\bm{x}(t),\bm{u}(t)) \, ,
$$
the dual of Problem P can be formulated as follows  \cite{hartl1995survey}.
\begin{prob}[Problem $P_\lambda$]
Determine ${\bm x}(t)$, ${\bm u}(t)$, ${\bm \lambda}(t)$, ${\bm \mu}(t)$ and ${\bm \nu}$ that for all $t \in [0,1]$ satisfy Equations \eqref{eq:dynamicconstraint}, \eqref{eq:equalityconstraint}, \eqref{eq:inequalityconstraint} and
\begin{align}  \label{eq:Plambda1}
& \bm{\mu}^\top (t) \bm{h}({\bm x}(t),{\bm u}(t)) = 0 \, , \quad   \bm{\mu} (t) \geq 0 \, , \\ \label{eq:Plambda2}
& \dot{\bm \lambda}^\top(t) + \mathcal{L}_x (\bm{x}(t),\bm{u}(t),\bm{\lambda}(t),\bm{\mu}(t)) =  0 \, ,
%\\ \nonumber & \quad  \dot{\bm \lambda}^\top(t) + F_{x}({\bm x}(t),{\bm u}(t)) +  {\bm \lambda}^\top(t) {\bm f}_{x} ({\bm x}(t),{\bm u}(t))
%\\ \nonumber & \qquad +{\bm \mu}^\top(t) {\bm h_{x}} ({\bm x}(t),{\bm u}(t)) = 0 \, , 
\\  \label{eq:Plambda3}
& \bm{\lambda}^\top(0) = -\bm{\nu}^\top \bm{e}_{x(0)}({\bm x}(0),{\bm x}(1)) - E_{x(0)}({\bm x}(0),{\bm x}(1)) \, , \\ \label{eq:Plambda4}
& \bm{\lambda}^\top(1) = \bm{\nu}^\top \bm{e}_{x(1)}({\bm x}(0),{\bm x}(1)) + E_{x(1)}({\bm x}(0),{\bm x}(1)) \, , \\ \label{eq:Plambda5}
& \mathcal{L}_u(\bm{x}(t),\bm{u}(t),\bm{\lambda}(t),\bm{\mu}(t)) = 0 \, .
% \bm{\lambda}^\top(t) \bm{f}_{u}({\bm x}(t),{\bm u}(t)) + F_{u}({\bm x}(t),{\bm u}(t)) + \bm{\mu}^\top(t) \bm{h}_{u}({\bm x}(t),{\bm u}(t)) = 0 \, .
\end{align}
\end{prob}
In the above problem, subscripts are used to denote partial derivatives, e.g. $F_{x}({\bm x},{\bm u}) = \frac{\partial}{\partial {\bm x}}F({\bm x},{\bm u})$.

The following assumptions are imposed onto Problem $P_\lambda$.
\begin{asm} \label{asm:functions}
$E,F,\bm{f},\bm{e}$ and $\bm{h}$ are continuously differentiable with respect to their arguments, and their gradients are Lipschitz continuous over the domain.
\end{asm}

\begin{asm} \label{asm:plambdaexists}
Solutions ${\bm x}^*(t)$, ${\bm u}^*(t)$, ${\bm \lambda}^*(t)$, ${\bm \mu}^*(t)$ and ${\bm \nu}^*$ of Problem $P_\lambda$ exist and satisfy ${\bm x}^*(t) \in \mathcal{C}^1_{n_x}$, $\bm{u}^*(t) \in \mathcal{C}^0_{n_u}$, ${\bm \lambda}^*(t) \in \mathcal{C}^1_{n_x}$ and $\bm{\mu}^*(t) \in \mathcal{C}^0_{n_h}$ in $[0,1]$.
\end{asm}

\begin{rem}
Notice that Problem $P_\lambda$ implicitly assumes the absence of pure state constraints in Problem $P$. If the inequality constraint in Equation \eqref{eq:inequalityconstraint} is independent of $\bm{u}(t)$, then the costate ${\bm \lambda}(t)$ must also satisfy the following jump condition \cite{hartl1995survey}:
$${\bm \lambda}(t_e^-) = {\bm \lambda}(t_e^+) + {\bm h}_{x(t_e)}^\top {\bm \eta} \, ,$$
% where $t_e$ is the entry or exit time into a constrained arc in which the inequality constraint is active, $t_e^-$ and $t_e^+$ denote the left-hand side and right-hand side limits of the trajectory, respectively, and ${\bm \eta}$ is a constant covector. For simplicity, the theoretical results that will be presented in Section \ref{sec:costateresults} do not consider the jump conditions above, i.e., the inequality constraints are dependent on $\bm{u}(t)$.
\end{rem}

\subsection{KKT conditions of Problem $P_N$}

Now we derive the relaxed KKT conditions associated with Problem $P_N$. Since the constraints \eqref{eq:discdynamicconstraint} and \eqref{eq:discinequalityconstraint} are imposed at the Bernstein coefficients, one multiplier is associated with each coefficient. Let $\blam_{j,N} \in \IR^{n_x}$ and $\bmu_{j,N} \in \IR^{n_h}$, $j=0,\ldots,N$, be the multipliers of the dynamic and inequality constraints at the $j$th coefficient, let $\bar{\bm{\nu}} \in \IR^{n_e}$ be the multiplier of the boundary conditions, and let
$\blam_N \in \IR^{n_x \times (N+1)}$ and $\bmu_N \in \IR^{n_h \times (N+1)}$ be defined as
$$\blam_N = [\blam_{0,N} \, , \, \ldots \, , \, \blam_{N,N}] , \qquad \bmu_N = [\bmu_{0,N} \, , \, \ldots \, , \, \bmu_{N,N}].$$
We also introduce the corresponding $N$th order Bernstein polynomials
\begin{equation} \label{eq:polappr_costates}
{\bm \lambda}_N(t) =   \sum_{j=0}^N \blam_{j,N} b_{j,N}(t) \, , \quad  {\bm \mu}_N(t) =   \sum_{j=0}^N \bmu_{j,N} b_{j,N}(t) \, .
\end{equation}
With the dynamic constraint written in its exact form $\sum_k\bx_{k,N}D_{k,j}-\bm{f}(\bx_{j,N},\bu_{j,N})=\bm{0}$, the Lagrangian of the equality-constrained version of Problem $P_N$ is
\begin{equation} \label{eq:LagrangianPN}
\begin{split}
\mathcal{L}_N & = E(\bx_{0,N},\bx_{N,N})+w \sum_{j=0}^N F(\bx_{j,N},\bu_{j,N})  + \sum_{j=0}^N \blam_{j,N}^\top \Big(\bm{f}(\bx_{j,N},\bu_{j,N}) - \sum_{k=0}^N \bx_{k,N}D_{k,j}\Big)
\\ & \qquad + \sum_{j=0}^N \bmu_{j,N}^\top {\bm h}(\bx_{j,N},\bu_{j,N}) + \bar{\bm{\nu}}^\top {\bm e}(\bx_{0,N},\bx_{N,N})\, .
\end{split}
\end{equation}
In what follows, a superscript $[j]$ denotes evaluation at the $j$th pair of coefficients, e.g., $F_x^{[j]} = F_x(\bx_{j,N},\bu_{j,N})$, and a superscript $[N]$ on a boundary term denotes evaluation at $(\bx_{0,N},\bx_{N,N})$; e.g., $E_{x(0)}^{[N]} = E_{x(0)}(\bx_{0,N},\bx_{N,N})$. Let $\delta_{k,j}$ be the Kronecker delta. Then, for $k=0,\ldots,N$,
\begin{align}
\frac{\partial \mathcal{L}_N}{\partial \bx_{k,N}} & = w F_x^{[k]} + \blam_{k,N}^\top\bm{f}_x^{[k]} + \bmu_{k,N}^\top\bm{h}_x^{[k]} - \sum_{j=0}^N D_{k,j}\blam_{j,N}^\top \nonumber \\ & \qquad + \delta_{k,0}\big(E_{x(0)}^{[N]}+\bar{\bm\nu}^\top\bm{e}_{x(0)}^{[N]}\big) + \delta_{k,N}\big(E_{x(1)}^{[N]}+\bar{\bm\nu}^\top\bm{e}_{x(1)}^{[N]}\big) \, , \label{eq:dLdx} \\
\frac{\partial \mathcal{L}_N}{\partial \bu_{k,N}} & = w F_u^{[k]} + \blam_{k,N}^\top\bm{f}_u^{[k]} + \bmu_{k,N}^\top\bm{h}_u^{[k]} \, . \label{eq:dLdu}
\end{align}
The structure of these derivatives becomes transparent in terms of the scaled multipliers
\begin{equation} \label{eq:deflammutilde}
\tlam_{j,N} = \frac{\blam_{j,N}}{w} \, , \qquad \tmu_{j,N} = \frac{\bmu_{j,N}}{w} \, , \qquad j=0,\ldots,N \, ,
\end{equation}
with corresponding Bernstein polynomials
$$\tilde{\bm \lambda}_N(t) = \sum_{j} \tlam_{j,N} b_{j,N}(t) = \frac{\bm{\lambda}_N(t)}{w} \, , \qquad \tilde{\bm \mu}_N(t) = \sum_{j} \tmu_{j,N} b_{j,N}(t) = \frac{\bm{\mu}_N(t)}{w} \, .$$ Let $\bm{\Delta}_i(\tlam_N) = N(\tlam_{i+1,N}-\tlam_{i,N})$, and let
$$
\tilde{\mathcal{L}}_x^{[k]} = F_x^{[k]} + \tlam_{k,N}^\top\bm{f}_x^{[k]} + \tmu_{k,N}^\top\bm{h}_x^{[k]} \, , \qquad \tilde{\mathcal{L}}_u^{[k]} = F_u^{[k]} + \tlam_{k,N}^\top\bm{f}_u^{[k]} + \tmu_{k,N}^\top\bm{h}_u^{[k]}
$$
be the derivatives of $\mathcal{L}$ evaluated at $(\bx_{k,N},\bu_{k,N},\tlam_{k,N},\tmu_{k,N})$.
\begin{lem}[Scaled KKT conditions] \label{lem:scaledKKT}
The derivatives \eqref{eq:dLdx} and \eqref{eq:dLdu} satisfy
\begin{align}
\frac{1}{w}\frac{\partial \mathcal{L}_N}{\partial \bx_{k,N}} & = \tilde{\mathcal{L}}_x^{[k]} + \frac{N-k+1}{N}\bm{\Delta}_{k-1}^\top(\tlam_N)+\frac{k+1}{N}\bm{\Delta}_{k}^\top(\tlam_N) \, , \qquad k=1,\ldots,N-1 \, , \label{eq:scaledinterior} \\
\frac{\partial \mathcal{L}_N}{\partial \bx_{0,N}} & = \tlam_{0,N}^\top + \bar{\bm\nu}^\top\bm{e}_{x(0)}^{[N]} + E_{x(0)}^{[N]} + w\,\tilde{\mathcal{L}}_x^{[0]} + \frac{w}{N}\bm{\Delta}_{0}^\top(\tlam_N) \, , \label{eq:scaledzero} \\
\frac{\partial \mathcal{L}_N}{\partial \bx_{N,N}} & = -\tlam_{N,N}^\top + \bar{\bm\nu}^\top\bm{e}_{x(1)}^{[N]} + E_{x(1)}^{[N]} + w\,\tilde{\mathcal{L}}_x^{[N]} + \frac{w}{N}\bm{\Delta}_{N-1}^\top(\tlam_N) \, , \label{eq:scaledN} \\
\frac{1}{w}\frac{\partial \mathcal{L}_N}{\partial \bu_{k,N}} & = \tilde{\mathcal{L}}_u^{[k]} \, , \qquad k=0,\ldots,N \, . \label{eq:scaledu}
\end{align}
\end{lem}
\textbf{Proof:} Substitute $\blam_{j,N}=w\tlam_{j,N}$ and $\bmu_{j,N}=w\tmu_{j,N}$ into Equations \eqref{eq:dLdx} and \eqref{eq:dLdu}, and use Lemma \ref{lem:Dproperties}(ii) and (iii), with $\bar{\bm{a}}_j=\tlam_{j,N}$, for the term $-w\sum_jD_{k,j}\tlam_{j,N}$.

Equation \eqref{eq:scaledinterior} is a discretization of the costate equation \eqref{eq:Plambda2} at the Bernstein coefficients, in which, by Lemma \ref{lem:Dproperties}, the term $\frac{N-k+1}{N}\bm{\Delta}_{k-1}+\frac{k+1}{N}\bm{\Delta}_{k}$ plays the role of $\dot{\bm{\lambda}}$. Equations \eqref{eq:scaledzero} and \eqref{eq:scaledN} are discretizations of the transversality conditions \eqref{eq:Plambda3} and \eqref{eq:Plambda4}, perturbed by terms of order $w$, and Equation \eqref{eq:scaledu} is the pointwise counterpart of \eqref{eq:Plambda5}.
With this setup, we define the following relaxed KKT system, obtained from \eqref{eq:LagrangianPN} by imposing each condition up to a tolerance.
\begin{prob}[Problem $P_{N\lambda}$]
Determine $\bx_N$, $\bu_N$, $\blam_N$, $\bmu_N$ and $\bar{\bm \nu}$ that satisfy the primal feasibility conditions, namely Equations \eqref{eq:discdynamicconstraint}, \eqref{eq:discequalityconstraint} and \eqref{eq:discinequalityconstraint}, the complementary slackness and dual feasibility conditions
\begin{equation} \label{eq:complementarity}
\left\vert \bmu_{k,N}^\top \bm{h}(\bx_{k,N},\bu_{k,N}) \right\vert  \leq w \delta_D^N \, , \qquad
\bmu_{k,N}  \geq  - w \delta_D^N  \bm{1} \, ,  \qquad  \forall k=0,\ldots,N \, ,
\end{equation}
and the stationarity conditions
\begin{equation} \label{eq:stationarity}
\begin{split}
& \left\Vert \frac{\partial \mathcal{L}_N}{\partial \bx_{k,N}} \right\Vert \leq w\delta_D^N \, , \quad k=1,\ldots,N-1 \, , \qquad \left\Vert \frac{\partial \mathcal{L}_N}{\partial \bx_{k,N}} \right\Vert \leq \delta_D^N \, , \quad k\in\{0,N\} \, , \\
& \left\Vert \frac{\partial \mathcal{L}_N}{\partial \bu_{k,N}} \right\Vert \leq w\delta_D^N \, , \quad k = 0,\ldots,N \, ,
\end{split}
\end{equation}
where $w=\frac{1}{N+1}$ and $\delta_D^N$ is a small positive number that depends on $N$ and  satisfies $\lim_{N \to \infty} \delta_D^N = 0$.
\end{prob}

\begin{rem} \label{rem:PNlambda}
If $\delta_P^N=\delta_D^N=0$, Problem $P_{N\lambda}$ is the KKT system of Problem $P_N^0$, the problem obtained from $P_N$ with the dynamic constraint \eqref{eq:discdynamicconstraint} written as an equality, and every KKT point of \(PN^0\) solves \(P_{N \lambda}\) for all \(\delta_P^N,\delta_D^N>0\). $\sum_k\bx_{k,N}D_{k,j}=\bm{f}(\bx_{j,N},\bu_{j,N})$. For $\delta_P^N>0$, Problem $P_{N\lambda}$ is not the KKT system of Problem $P_N$, whose dynamic constraint is then an inequality on a norm. Similarly to the relaxation bound $\delta_P^N$, the bound $\delta_D^N$ is introduced to guarantee that Problem $P_{N\lambda}$ is feasible (see Theorem \ref{thm:existence_dual}). The relaxation in \eqref{eq:discdynamicconstraint}  cannot be removed in general, since \(D_N\) is nilpotent and \eqref{eq:discdynamicconstraint} with \(\delta_P^N=0\) can be infeasible. 
\end{rem}

\begin{rem}[Scaling of the tolerances] \label{rem:scaling}
By Lemma \ref{lem:scaledKKT}, the derivatives $\frac{\partial \mathcal{L}_N}{\partial \bx_{k,N}}$, $k=1,\ldots,N-1$, and $\frac{\partial \mathcal{L}_N}{\partial \bu_{k,N}}$, $k=0,\ldots,N$, are $w$ times quantities that are bounded whenever the scaled multipliers and the slopes $\bm{\Delta}_i(\tlam_N)$ are bounded. The corresponding tolerances in Equations \eqref{eq:complementarity} and \eqref{eq:stationarity} are therefore scaled by $w$. Without this scaling, these conditions would be satisfied by any bounded sequence of scaled multipliers for sufficiently large $N$, and would carry no information on Equations \eqref{eq:Plambda1}, \eqref{eq:Plambda2} and \eqref{eq:Plambda5} in the limit. In contrast, the derivatives with respect to $\bx_{0,N}$ and $\bx_{N,N}$ contain terms of order one, see Equations \eqref{eq:scaledzero} and \eqref{eq:scaledN}, and their tolerances are not scaled.
\end{rem}

The solution of Problem $P_{N\lambda}$ presents a set of Bernstein coefficients $\bx^*_N$, $\bu^*_N$, $\blam^*_N$, $\bmu^*_N$ (which determine the Bernstein polynomials $\bm{x}_N^*(t)$, $\bm{u}_N^*(t)$, $\bm{\lambda}_N^*(t)$ and $\bm{\mu}_N^*(t)$) and a vector $\bar{\bm \nu}^*$.

\section{Feasibility and consistency of Problem $P_{N \lambda}$} \label{sec:costateresults}

The objective of this section is to investigate the ability of the solutions of Problem $P_{N\lambda}$ to approximate the solutions of Problem $P_\lambda$.
In what follows, we first show the existence of a solution for Problem $P_{N\lambda}$ (feasibility).
Second, we investigate the convergence properties of the solutions of Problem $P_{N\lambda}$ as $N \to \infty$ (consistency). Third, by combining these two results, we formulate the \emph{covector mapping theorem} for Bernstein approximations, which relates the solutions of Problem $P_{N\lambda}$ to the solutions of Problem $P_\lambda$.
The main results of this section are reported in the three theorems below and summarized in Figure \ref{fig:covectormapping}. Throughout this section, $t_j=\frac{j}{N}$ and $w=\frac{1}{N+1}$.

\begin{thm}[Feasibility] \label{thm:existence_dual}
Let $(\bm{x}(t),\bm{u}(t),\bm{\lambda}(t),\bm{\mu}(t),\bm{\nu})$ be a solution of Problem $P_\lambda$ satisfying Assumption \ref{asm:plambdaexists}, and let
\begin{align}
\delta_P^N & = C_P W_{x^\prime}(N^{-1}) \, , \label{eq:defdeltaPN} \\
\delta_D^N & = C_D \max \left\{ W_{\lambda^\prime}(N^{-1}) \, , \, N^{-1} \right\} \, , \label{eq:defdeltaDN}
\end{align}
where $C_P$ and $C_D$ are positive constants independent of $N$, and $W_{x^\prime}(\cdot)$ and $W_{\lambda^\prime}(\cdot)$ are the moduli of continuity of $\dot{\bm{x}}(t)$ and $\dot{\bm{\lambda}}(t)$, respectively.
Then, for arbitrary order of approximation $N \in \mathbb{Z}^+$, Problem $P_{N\lambda}$ admits the solution
\begin{equation} \label{eq:defxu}
\bx_{j,N} = \bm{x}(t_j) \, , \quad \bu_{j,N} = \bm{u}(t_j) \, ,
\end{equation}
\begin{equation} \label{eq:deflammu}
{\blam}_{j,N} = w \bm{\lambda}(t_j) \, , \quad  {\bmu}_{j,N} = w \bm{\mu}(t_j) \, , \quad \bar{\bm{\nu}} = \bm{\nu} \, ,
\end{equation}
$j = 0 , \ldots , N$. Furthermore, the corresponding Bernstein polynomials satisfy, for all $t\in[0,1]$,
\begin{equation} \label{eq:boundedcostates}
\begin{split}
& ||{\bm{x}}_N(t)-\bm{x}(t)|| \leq C_x W_{x}(N^{-\frac{1}{2}}) \, , \qquad ||{\bm{u}}_N(t)-\bm{u}(t)|| \leq C_u W_{u}(N^{-\frac{1}{2}}) \, , \\
& \left\Vert\frac{\bm{\lambda}_N(t)}{w}-\bm{\lambda}(t)\right\Vert \leq C_\lambda W_{\lambda}(N^{-\frac{1}{2}}) \, , \qquad \left\Vert\frac{\bm{\mu}_N(t)}{w}-\bm{\mu}(t)\right\Vert \leq C_\mu W_{\mu}(N^{-\frac{1}{2}}) \, ,
\end{split}
\end{equation}
where $C_x$, $C_u$, $C_\lambda$ and $C_\mu$ are positive constants independent of $N$, and $W_{x}(\cdot)$, $W_{u}(\cdot)$, $W_{\lambda}(\cdot)$ and $W_{\mu}(\cdot)$ are the moduli of continuity of $\bm{x}(t)$, $\bm{u}(t)$, $\bm{\lambda}(t)$ and $\bm{\mu}(t)$, respectively.
\end{thm}

\textbf{Proof:} The proof is constructive. We show that the point defined by Equations \eqref{eq:defxu} and \eqref{eq:deflammu} satisfies each condition of Problem $P_{N\lambda}$. By Equation \eqref{eq:deflammutilde}, the scaled multipliers are the samples of the multipliers of Problem $P_\lambda$, i.e., $\tlam_{j,N}=\bm{\lambda}(t_j)$ and $\tmu_{j,N}=\bm{\mu}(t_j)$. Hence, all quantities with superscript $[j]$ in Lemma \ref{lem:scaledKKT} are evaluated along the solution of Problem $P_\lambda$ at $t_j$; e.g., $\tilde{\mathcal{L}}_x^{[j]}=\mathcal{L}_x(\bm{x}(t_j),\bm{u}(t_j),\bm{\lambda}(t_j),\bm{\mu}(t_j))$. Moreover, the quantities with superscript $[N]$ are evaluated at $(\bm{x}(0),\bm{x}(1))$, since $\bx_{0,N}=\bm{x}(0)$ and $\bx_{N,N}=\bm{x}(1)$.

\emph{Primal feasibility.} Equations \eqref{eq:discequalityconstraint} and \eqref{eq:discinequalityconstraint} are satisfied exactly, since $\bm{e}(\bx_{0,N},\bx_{N,N}) = \bm{e}(\bm{x}(0),\bm{x}(1))=\bm{0}$ and
$$\bm{h}(\bx_{j,N},\bu_{j,N}) = \bm{h}(\bm{x}(t_j),\bm{u}(t_j))\leq\bm{0} \, .$$ For the dynamic constraint, since $\bm{f}(\bm{x}(t_j),\bm{u}(t_j))=\dot{\bm{x}}(t_j)$, Equation \eqref{eq:Dright} in Lemma \ref{lem:Dproperties} gives
\begin{equation*}
\Big\Vert \sum_{k=0}^N \bx_{k,N}D_{k,j} - \bm{f}(\bx_{j,N},\bu_{j,N})\Big\Vert = \Big\Vert \sum_{k=0}^N \bm{x}(t_k)D_{k,j} - \dot{\bm{x}}(t_j)\Big\Vert \leq \sqrt{n_x}W_{x^\prime}(N^{-1}) \, .
\end{equation*}
Thus, Equation \eqref{eq:discdynamicconstraint} holds with $\delta_P^N$ given by Equation \eqref{eq:defdeltaPN} and $C_P\geq\sqrt{n_x}$.

\emph{Complementary slackness and dual feasibility.} By Equation \eqref{eq:Plambda1},
$$\bmu_{k,N}^\top \bm{h}(\bx_{k,N},\bu_{k,N}) = w\bm{\mu}^\top(t_k)\bm{h}(\bm{x}(t_k),\bm{u}(t_k))=0$$
and $\bmu_{k,N}=w\bm{\mu}(t_k)\geq\bm{0}$. Hence Equation \eqref{eq:complementarity} holds for any $\delta_D^N\geq 0$.

\emph{Stationarity with respect to $\bu_{k,N}$.} By Equations \eqref{eq:scaledu} and \eqref{eq:Plambda5}, $\frac{\partial \mathcal{L}_N}{\partial \bu_{k,N}} = w\,\tilde{\mathcal{L}}_u^{[k]}=\bm{0}$, $k=0,\ldots,N$, since $\tilde{\mathcal{L}}_u^{[k]}$ is the left hand side of Equation \eqref{eq:Plambda5} at $t=t_k$.

\emph{Stationarity with respect to $\bx_{k,N}$, $k=1,\ldots,N-1$.} By Equation \eqref{eq:Plambda2}, $\tilde{\mathcal{L}}_x^{[k]} = -\dot{\bm{\lambda}}^\top(t_k)$. Then, by Equation \eqref{eq:scaledinterior}, Lemma \ref{lem:Dproperties}(ii) and Equation \eqref{eq:Dleft} with $\bm{a}(t)=\bm{\lambda}(t)$,
$$
\frac{1}{w}\left\Vert\frac{\partial \mathcal{L}_N}{\partial \bx_{k,N}}\right\Vert = \Big\Vert -\sum_{j=0}^N D_{k,j}\bm{\lambda}(t_j) - \dot{\bm{\lambda}}(t_k) \Big\Vert \leq \sqrt{n_x}W_{\lambda^\prime}(N^{-1}) + 2\sqrt{n_x}||\dot{\bm{\lambda}}||_\infty N^{-1} \, .
$$

\emph{Stationarity with respect to $\bx_{0,N}$ and $\bx_{N,N}$.} By Equation \eqref{eq:Plambda3},
$$
\tlam_{0,N}^\top + \bar{\bm\nu}^\top\bm{e}_{x(0)}^{[N]} + E_{x(0)}^{[N]} = \bm{\lambda}^\top(0) + \bm{\nu}^\top\bm{e}_{x(0)}(\bm{x}(0),\bm{x}(1)) + E_{x(0)}(\bm{x}(0),\bm{x}(1)) = \bm{0} \, .
$$
Then Equation \eqref{eq:scaledzero}, together with $\tilde{\mathcal{L}}_x^{[0]}=-\dot{\bm{\lambda}}^\top(0)$, $||\bm{\Delta}_0(\tlam_N)||\leq\sqrt{n_x}||\dot{\bm{\lambda}}||_\infty$ (Lemma \ref{lem:Dproperties}) and $w\leq N^{-1}$, gives
$$
\left\Vert\frac{\partial \mathcal{L}_N}{\partial \bx_{0,N}}\right\Vert = \left\Vert -w\dot{\bm{\lambda}}^\top(0) + \frac{w}{N}\bm{\Delta}_0^\top(\tlam_N) \right\Vert \leq (1+\sqrt{n_x})||\dot{\bm{\lambda}}||_\infty N^{-1} \, .
$$
Using Equations \eqref{eq:Plambda4} and \eqref{eq:scaledN}, the same bound is obtained for $\bx_{N,N}$.

Combining the bounds above, the conditions \eqref{eq:complementarity} and \eqref{eq:stationarity} are satisfied with $\delta_D^N$ given by Equation \eqref{eq:defdeltaDN}, provided that $C_D \geq \sqrt{n_x} + (1+2\sqrt{n_x})||\dot{\bm{\lambda}}||_\infty$, which is independent of $N$. Since $\dot{\bm{x}}(t)$ and $\dot{\bm{\lambda}}(t)$ are continuous, $\delta_P^N\to0$ and $\delta_D^N\to0$ as $N\to\infty$.
Finally, by Equations \eqref{eq:defxu} and \eqref{eq:deflammu}, $\bm{x}_N(t)$, $\bm{u}_N(t)$, $\frac{\bm{\lambda}_N(t)}{w}=\tilde{\bm{\lambda}}_N(t)$ and $\frac{\bm{\mu}_N(t)}{w}=\tilde{\bm{\mu}}_N(t)$ are the Bernstein approximations of $\bm{x}(t)$, $\bm{u}(t)$, $\bm{\lambda}(t)$ and $\bm{\mu}(t)$, and Equation \eqref{eq:boundedcostates} follows from Lemma \ref{lem:bernsteinapprox}.

\begin{cor} \label{cor:c2_dual}
	If the solution of Problem $P_\lambda$ in Theorem \ref{thm:existence_dual} satisfies ${\bm x}(t) \in \mathcal{C}^2_{n_x}$ and ${\bm \lambda}(t) \in \mathcal{C}^2_{n_x}$ on $[0,1]$, then Theorem \ref{thm:existence_dual} holds with
	$ \delta_{P}^N = C_P N^{-1} $ and $ \delta_{D}^N = C_D N^{-1}$,
	where $C_P$ and $C_D$ are positive constants independent of $N$. If, in addition, $\bm{u}(t)$ and $\bm{\mu}(t)$ are of class $\mathcal{C}^2$, then the right hand sides of Equation \eqref{eq:boundedcostates} can be replaced by $C_xN^{-1}$, $C_uN^{-1}$, $C_\lambda N^{-1}$ and $C_\mu N^{-1}$.
\end{cor}

\textbf{Proof:} Under the stated regularity, $W_{x^\prime}(\delta)\leq ||\ddot{\bm{x}}||_\infty\delta$ and $W_{\lambda^\prime}(\delta)\leq ||\ddot{\bm{\lambda}}||_\infty\delta$, and the first statement follows from Equations \eqref{eq:defdeltaPN} and \eqref{eq:defdeltaDN}. The second statement follows by applying Lemma \ref{lem:bernsteinapproxC2} instead of Lemma \ref{lem:bernsteinapprox} in the last step of the proof of Theorem \ref{thm:existence_dual}.

\begin{rem}
For arbitrarily small scalars $\epsilon_P>0$ and $\epsilon_D>0$ there exists $N_1$ such that for all $N \geq N_1$, we have $\delta_P^N \leq \epsilon_P$ and $\delta_D^N \leq \epsilon_D$; i.e., the relaxation bounds in Problem $P_{N \lambda}$ can be made arbitrarily small by choosing sufficiently large $N$. Notice that the constraint residuals at the Bernstein coefficients are of order $W(N^{-1})$, whereas the uniform distance between the Bernstein polynomials and the functions they approximate in Equation \eqref{eq:boundedcostates} is of order $W(N^{-\frac{1}{2}})$, or $N^{-1}$ under the assumptions of Corollary \ref{cor:c2_dual}.
\end{rem}

We now turn to the consistency of Problem $P_{N\lambda}$. Following the approach of \cite{cichella2020optimal}, we formulate the required convergence assumption in terms of control polygons. Let
$$
\{(\bx_N^*,\bu_N^*,\blam_N^*,\bmu_N^*,\bar{\bm{\nu}}_N^*)\}_{N=N_1}^\infty
$$
be a sequence of solutions of Problem $P_{N\lambda}$, and let
\begin{equation} \label{eq:constistencymappingone}
\tlam_{j,N}^* = \frac{{\blam}_{j,N}^*}{w} \, , \quad \tmu_{j,N}^* = \frac{{\bmu}_{j,N}^*}{w} \, , \qquad j=0,\ldots,N \, .
\end{equation}
Let $\bm{X}_N^*(t)$, $\bm{U}_N^*(t)$, $\tilde{\bm{\Lambda}}_N^*(t)$ and $\tilde{\bm{M}}_N^*(t)$, $t\in[0,1]$, be the control polygons defined by $\bx_N^*$, $\bu_N^*$, $\tlam_N^*$ and $\tmu_N^*$, respectively (see Equation \eqref{eq:controlpolygon}).

\begin{asm} \label{asm:primalconv}
There exist $\bm{x}^\infty(t)\in\mathcal{C}^1_{n_x}$ and $\bm{u}^\infty(t)\in\mathcal{C}^0_{n_u}$ on $[0,1]$ such that
\begin{equation*}
\begin{split}
& \lim_{N\to\infty}\max_{t\in[0,1]}\Big( ||\bm{X}_N^*(t)-\bm{x}^\infty(t)|| + ||\bm{U}_N^*(t)-\bm{u}^\infty(t)|| \Big) = 0 \, , \\
& \lim_{N\to\infty}\max_{i\in\{0,\ldots,N-1\}} \Big\Vert \bm{\Delta}_i(\bx_N^*) - \dot{\bm{x}}^\infty(t_i)\Big\Vert = 0 \, ,
\end{split}
\end{equation*}
where $\bm{\Delta}_i(\bx_N^*)=N(\bx_{i+1,N}^*-\bx_{i,N}^*)$.
\end{asm}

\begin{asm} \label{asm:dualconv}
There exist $\bm{\lambda}^\infty(t)\in\mathcal{C}^1_{n_x}$ and $\bm{\mu}^\infty(t)\in\mathcal{C}^0_{n_h}$ on $[0,1]$ and $\bm{\nu}^\infty\in\IR^{n_e}$ such that $\lim_{N\to\infty}\bar{\bm{\nu}}_N^*=\bm{\nu}^\infty$, and
\begin{equation*}
\begin{split}
& \lim_{N\to\infty}\max_{t\in[0,1]}\Big( ||\tilde{\bm{\Lambda}}_N^*(t)-\bm{\lambda}^\infty(t)|| + ||\tilde{\bm{M}}_N^*(t)-\bm{\mu}^\infty(t)|| \Big) = 0 \, , \\
& \lim_{N\to\infty}\max_{i\in\{0,\ldots,N-1\}} \Big\Vert \bm{\Delta}_i(\tlam_N^*) - \dot{\bm{\lambda}}^\infty(t_i)\Big\Vert = 0 \, ,
\end{split}
\end{equation*}
where $\bm{\Delta}_i(\tlam_N^*)=N(\tlam_{i+1,N}^*-\tlam_{i,N}^*)$.
\end{asm}

Assumption \ref{asm:primalconv} is stated on the control polygons and on their slopes, as in the consistency analysis of \cite{cichella2020optimal}, and Assumption \ref{asm:dualconv} is its counterpart for the scaled multipliers. In simple words, these assumptions require that, as $N$ increases, the solutions of Problem $P_{N\lambda}$ settle on a limit, and that the slopes of the control polygons of the state and of the scaled costate converge to the derivatives of the limit. Since the control polygons interpolate the coefficients, i.e., $\bm{X}_N^*(t_k)=\bx_{k,N}^*$, the assumptions imply
\begin{equation} \label{eq:coeffconv}
\begin{split}
\eta_N := \max_{k\in\{0,\ldots,N\}}\Big( & ||\bx_{k,N}^*-\bm{x}^\infty(t_k)|| + ||\bu_{k,N}^*-\bm{u}^\infty(t_k)|| \\
& + ||\tlam_{k,N}^*-\bm{\lambda}^\infty(t_k)|| + ||\tmu_{k,N}^*-\bm{\mu}^\infty(t_k)|| \Big) \to 0 \, .
\end{split}
\end{equation}

\begin{thm}[Consistency] \label{thm:consistency_dual}
Let $\{(\bx_N^*,\bu_N^*,\blam_N^*,\bmu_N^*,\bar{\bm{\nu}}_N^*)\}_{N=N_1}^\infty$ be a sequence of solutions of Problem $P_{N\lambda}$ satisfying Assumptions \ref{asm:primalconv} and \ref{asm:dualconv}, with $\lim_{N\to\infty}\delta_P^N=\lim_{N\to\infty}\delta_D^N=0$, and let Assumption \ref{asm:functions} hold. Then
$$({\bm x}^\infty(t),{\bm u}^\infty(t), {\bm \lambda}^\infty(t),{\bm \mu}^\infty(t),{\bm \nu}^\infty)$$
is a solution of Problem $P_\lambda$. Moreover, the Bernstein polynomials ${\bm x}_N^*(t)$, ${\bm u}_N^*(t)$, $\tilde{\bm \lambda}_N^*(t)=\frac{\bm{\lambda}_N^*(t)}{w}$ and $\tilde{\bm \mu}_N^*(t)=\frac{\bm{\mu}_N^*(t)}{w}$ converge uniformly on $[0,1]$ to this solution.
\end{thm}

\textbf{Proof:} The objective is to show that the limit satisfies Equations \eqref{eq:dynamicconstraint}--\eqref{eq:inequalityconstraint} and \eqref{eq:Plambda1}--\eqref{eq:Plambda5}. Throughout the proof, for a given $t'\in[0,1]$, $k_N=\lfloor Nt'\rfloor$, so that $\lim_{N\to\infty}t_{k_N}=t'$. Let
$$
\varepsilon_N := \max_{i\in\{0,\ldots,N-1\}}\Big(\Big\Vert \bm{\Delta}_i(\bx_N^*) - \dot{\bm{x}}^\infty(t_i)\Big\Vert + \Big\Vert \bm{\Delta}_i(\tlam_N^*) - \dot{\bm{\lambda}}^\infty(t_i)\Big\Vert\Big) \, ,
$$
which converges to zero by Assumptions \ref{asm:primalconv} and \ref{asm:dualconv}. By Equation \eqref{eq:coeffconv}, the coefficients $\bx_{k,N}^*$, $\bu_{k,N}^*$, $\tlam_{k,N}^*$ and $\tmu_{k,N}^*$, as well as the slopes $\bm{\Delta}_i(\bx_N^*)$ and $\bm{\Delta}_i(\tlam_N^*)$, are bounded uniformly in $N$ and $k$. Hence, there is a compact set $K$ containing all the arguments at which the problem data are evaluated, on which, by Assumption \ref{asm:functions}, $E$, $F$, $\bm{f}$, $\bm{e}$, $\bm{h}$ and their gradients are bounded and Lipschitz continuous. In particular, since $\tilde{\mathcal{L}}_x^{[k]}$ and $\tilde{\mathcal{L}}_u^{[k]}$ are affine in $(\tlam_{k,N}^*,\tmu_{k,N}^*)$ with coefficients that are Lipschitz in $(\bx_{k,N}^*,\bu_{k,N}^*)$, there exists $L>0$, independent of $N$, such that
\begin{equation} \label{eq:Lipschitzlimit}
\Big\Vert \tilde{\mathcal{L}}_x^{[k]} - \mathcal{L}_x^\infty(t_k) \Big\Vert + \Big\Vert \tilde{\mathcal{L}}_u^{[k]} - \mathcal{L}_u^\infty(t_k) \Big\Vert \leq L\eta_N \, , \qquad k=0,\ldots,N \, ,
\end{equation}
where
$$\mathcal{L}_x^\infty(t)=\mathcal{L}_x(\bm{x}^\infty(t),\bm{u}^\infty(t),\bm{\lambda}^\infty(t),\bm{\mu}^\infty(t)) \, ,$$
and similarly for $\mathcal{L}_u^\infty(t)$ and $\bm{f}^\infty(t)=\bm{f}(\bm{x}^\infty(t),\bm{u}^\infty(t))$.

\emph{Step 1: primal feasibility.} By Lemma \ref{lem:Dproperties}(i), the weights in $\sum_i\bx_{i,N}^*D_{i,k}$ are nonnegative, sum to one and multiply the slopes on the segments adjacent to the $k$th node, hence
\begin{equation} \label{eq:discretederivativelimit}
\Big\Vert \sum_{i=0}^N \bx_{i,N}^*D_{i,k} - \dot{\bm{x}}^\infty(t_k)\Big\Vert \leq \varepsilon_N + \sqrt{n_x}W_{\dot{x}^\infty}(N^{-1}) \, , \qquad k=0,\ldots,N \, .
\end{equation}
Then, by the triangle inequality, Equation \eqref{eq:discdynamicconstraint} and the Lipschitz continuity of $\bm{f}$,
$$
||\dot{\bm{x}}^\infty(t_{k_N}) - \bm{f}^\infty(t_{k_N})|| \leq \varepsilon_N + \sqrt{n_x}W_{\dot{x}^\infty}(N^{-1}) + \delta_P^N + L_f\eta_N \, ,
$$
where $L_f$ is the Lipschitz constant of $\bm{f}$ on $K$. The right hand side vanishes as $N\to\infty$, while, by continuity, the left hand side converges to $||\dot{\bm{x}}^\infty(t') - \bm{f}^\infty(t')||$. Since $t'$ is arbitrary, Equation \eqref{eq:dynamicconstraint} holds. Equations \eqref{eq:equalityconstraint} and \eqref{eq:inequalityconstraint} follow from
\begin{gather*}
\bm{e}(\bm{x}^\infty(0),\bm{x}^\infty(1)) = \lim_{N\to\infty}\bm{e}(\bx_{0,N}^*,\bx_{N,N}^*) = \bm{0} \, , \\ \bm{h}(\bm{x}^\infty(t'),\bm{u}^\infty(t')) = \lim_{N\to\infty}\bm{h}(\bx_{k_N,N}^*,\bu_{k_N,N}^*)\leq\bm{0} \, ,
\end{gather*}
where we used Equations \eqref{eq:coeffconv}, \eqref{eq:discequalityconstraint} and \eqref{eq:discinequalityconstraint} and the continuity of $\bm{e}$, $\bm{h}$, $\bm{x}^\infty$ and $\bm{u}^\infty$.

\emph{Step 2: Equation \eqref{eq:Plambda1}.} Dividing Equation \eqref{eq:complementarity} by $w$ gives $|\tmu_{k,N}^{*\top}\bm{h}(\bx_{k,N}^*,\bu_{k,N}^*)|\leq\delta_D^N$ and $\tmu_{k,N}^*\geq-\delta_D^N\bm{1}$. Evaluating at $k=k_N$ and letting $N\to\infty$, Equation \eqref{eq:coeffconv} and the continuity of $\bm{h}$, $\bm{x}^\infty$, $\bm{u}^\infty$ and $\bm{\mu}^\infty$ yield $\bm{\mu}^{\infty\top}(t')\bm{h}(\bm{x}^\infty(t'),\bm{u}^\infty(t'))=0$ and $\bm{\mu}^\infty(t')\geq\bm{0}$.

\emph{Step 3: Equation \eqref{eq:Plambda2}.} Let $t'\in(0,1)$, so that $1\leq k_N\leq N-1$ for sufficiently large $N$. Let $\bm{g}_k := \frac{N-k+1}{N}\bm{\Delta}_{k-1}(\tlam_N^*)+\frac{k+1}{N}\bm{\Delta}_{k}(\tlam_N^*) = \frac{N+2}{N}\big(\alpha_k\bm{\Delta}_{k-1}(\tlam_N^*)+(1-\alpha_k)\bm{\Delta}_{k}(\tlam_N^*)\big)$, with $\alpha_k=\frac{N-k+1}{N+2}\in(0,1)$. Then
$$
||\bm{g}_k - \dot{\bm{\lambda}}^\infty(t_k)|| \leq \frac{N+2}{N}\Big(\varepsilon_N + \sqrt{n_x}W_{\dot{\lambda}^\infty}(N^{-1})\Big) + \frac{2}{N}||\dot{\bm{\lambda}}^\infty||_\infty \, .
$$
By Equations \eqref{eq:scaledinterior} and \eqref{eq:stationarity}, $||\tilde{\mathcal{L}}_x^{[k]}+\bm{g}_k^\top||\leq\delta_D^N$ for $k=1,\ldots,N-1$. Hence, using Equation \eqref{eq:Lipschitzlimit},
$$
||\dot{\bm{\lambda}}^{\infty\top}(t_{k_N}) + \mathcal{L}_x^\infty(t_{k_N})|| \leq \delta_D^N + L\eta_N + \frac{N+2}{N}\Big(\varepsilon_N + \sqrt{n_x}W_{\dot{\lambda}^\infty}(N^{-1})\Big) + \frac{2}{N}||\dot{\bm{\lambda}}^\infty||_\infty \, .
$$
The right hand side vanishes as $N\to\infty$, and the left hand side converges to $||\dot{\bm{\lambda}}^{\infty\top}(t') + \mathcal{L}_x^\infty(t')||$ by continuity. Thus Equation \eqref{eq:Plambda2} holds on $(0,1)$, and, by continuity of both terms, on $[0,1]$.

\emph{Step 4: Equation \eqref{eq:Plambda5}.} By Equations \eqref{eq:scaledu}, \eqref{eq:stationarity} and \eqref{eq:Lipschitzlimit}, $||\mathcal{L}_u^\infty(t_{k_N})||\leq\delta_D^N+L\eta_N$ for all $t'\in[0,1]$. Letting $N\to\infty$ gives $\mathcal{L}_u^\infty(t')=\bm{0}$.

\emph{Step 5: Equations \eqref{eq:Plambda3} and \eqref{eq:Plambda4}.} By Equations \eqref{eq:scaledzero} and \eqref{eq:stationarity},
$$
\left\Vert \tlam_{0,N}^{*\top} + \bar{\bm{\nu}}_N^{*\top}\bm{e}_{x(0)}^{[N]} + E_{x(0)}^{[N]} \right\Vert \leq \delta_D^N + w\Big(||\tilde{\mathcal{L}}_x^{[0]}|| + \frac{1}{N}||\bm{\Delta}_0(\tlam_N^*)||\Big) \, ,
$$
where the superscript $[N]$ denotes evaluation at $(\bx_{0,N}^*,\bx_{N,N}^*)$. The terms in parentheses are bounded uniformly in $N$, and $w\to0$; hence the right hand side vanishes. By Equation \eqref{eq:coeffconv}, Assumption \ref{asm:dualconv}, and the continuity of $E_{x(0)}$ and $\bm{e}_{x(0)}$, the left hand side converges to $||\bm{\lambda}^{\infty\top}(0)+\bm{\nu}^{\infty\top}\bm{e}_{x(0)}(\bm{x}^\infty(0),\bm{x}^\infty(1))+E_{x(0)}(\bm{x}^\infty(0),\bm{x}^\infty(1))||$, which proves Equation \eqref{eq:Plambda3}. Equation \eqref{eq:Plambda4} follows identically from Equation \eqref{eq:scaledN}.

\emph{Convergence of the polynomials.} Since the basis polynomials are nonnegative and sum to one,
$
||\tilde{\bm{\lambda}}_N^*(t)-\bm{\lambda}^\infty(t)|| \leq \eta_N + ||\sum_j\bm{\lambda}^\infty(t_j)b_{j,N}(t)-\bm{\lambda}^\infty(t)||
$,
where the last term is the error of the Bernstein approximation of $\bm{\lambda}^\infty$, which vanishes uniformly by Lemma \ref{lem:bernsteinapprox}. The same argument applies to $\bm{x}_N^*$, $\bm{u}_N^*$ and $\tilde{\bm{\mu}}_N^*$ (see also Lemma \ref{lem:controlpolygon}). This completes the proof of Theorem \ref{thm:consistency_dual}.

\begin{rem}
Theorem \ref{thm:consistency_dual} does not require the primal sequence to be optimal for Problem $P_N$; it shows that the limit is an extremal of Problem $P$, i.e., a solution of the first order conditions in Problem $P_\lambda$. If, in addition, $(\bx_N^*,\bu_N^*)$ is an optimal solution of Problem $P_N$ for each $N$, then, under the assumptions of \cite{cichella2020optimal}, the limit $(\bm{x}^\infty(t),\bm{u}^\infty(t))$ is an optimal solution of Problem $P$ \cite{cichella2020optimal}.
\end{rem}

\begin{thm}[Covector Mapping Theorem] \label{thm:covector}
Let Assumption \ref{asm:functions} hold. Consider the \emph{covector mapping}
\begin{equation} \label{eq:mapping}
\begin{split}
& \bx_{j,N}\mapsto\bm{x}(t_j) \, , \quad \bu_{j,N}\mapsto\bm{u}(t_j) \, , \\
& (N+1)\,\blam_{j,N}\mapsto\bm{\lambda}(t_j) \, , \quad (N+1)\,\bmu_{j,N}\mapsto\bm{\mu}(t_j) \, , \quad \bar{\bm \nu} \mapsto \bm{\nu} \, ,
\end{split}
\end{equation}
$j=0,\ldots,N$, between the variables of Problem $P_{N\lambda}$ and those of Problem $P_\lambda$. Then:
\begin{enumerate}
\item[(i)] for any solution of Problem $P_\lambda$ satisfying Assumption \ref{asm:plambdaexists}, the inverse of the mapping \eqref{eq:mapping} defines a solution of Problem $P_{N\lambda}$, with $\delta_P^N$ and $\delta_D^N$ given by Equations \eqref{eq:defdeltaPN} and \eqref{eq:defdeltaDN}, for every $N\in\mathbb{Z}^+$;
\item[(ii)] for any sequence of solutions of Problem $P_{N\lambda}$ satisfying Assumptions \ref{asm:primalconv} and \ref{asm:dualconv}, the images under the mapping \eqref{eq:mapping} converge, as $N \to \infty$, to a solution of Problem $P_\lambda$; in particular, $\frac{\bm{\lambda}_N^*(t)}{w}$, $\frac{\bm{\mu}_N^*(t)}{w}$ and $\bar{\bm{\nu}}_N^*$ converge to its costate and multipliers.
\end{enumerate}
\end{thm}

\textbf{Proof:} Statement (i) is Theorem \ref{thm:existence_dual}, since the inverse of the mapping \eqref{eq:mapping} is given by Equations \eqref{eq:defxu} and \eqref{eq:deflammu}. Statement (ii) is Theorem \ref{thm:consistency_dual}, since $(N+1)\blam_{j,N}^*=\tlam_{j,N}^*$ and $(N+1)\bmu_{j,N}^*=\tmu_{j,N}^*$.

\begin{rem}  \label{rem:hamiltonian}
Define the Hamiltonian approximation
\[
\mathcal{H}_N(t)= F({\bm x_N^*}(t),{\bm u_N^*}(t)) + \frac{\bm{\lambda}_{N}^{*\top}(t)}{w} \bm{f}(\bm{x}_N^*(t),\bm{u}_N^*(t)) \, .
\]
Under the assumptions of Theorem \ref{thm:consistency_dual}, the uniform convergence of the polynomials and the continuity of $F$ and $\bm{f}$ imply
$
\lim_{N \rightarrow \infty} \max_{t\in[0,1]}|\mathcal{H}_N(t) - \mathcal{H}(\bm{x}^\infty(t),\bm{u}^\infty(t),\bm{\lambda}^\infty(t))| = 0
$.
\end{rem}
%==============================================================================

\tikzstyle{block} = [draw, fill=blue!20, rectangle, 
    minimum height=3em, minimum width=4em]
\tikzstyle{sum} = [draw, fill=blue!20, circle]
\tikzstyle{input} = [coordinate]

\begin{figure} 
\centering
%\resizebox{\columnwidth}{!}{%
\begin{tikzpicture}
    \node (B) [block, minimum width=3.0cm, minimum height=1cm] {Problem $P$};
    \node (C) [block, right of=B,xshift=6cm,minimum width=3.0cm, minimum height=1cm] {Problem $P_N$};
    \draw [-latex] ([yshift=+0.30cm]B.east) -- node[pos=0.5,above] {\scriptsize{Bernstein approximation}} ([yshift=+0.30cm]C.west);
    \draw [-latex] ([yshift=-0.30cm]C.west) -- node[pos=0.5,above] {\scriptsize{Convergence}} ([yshift=-0.30cm]B.east);
    \node [block, below of=C, yshift=-1.8cm,minimum width=3.0cm, minimum height=1cm] (Clam) {Problem $P_{N \lambda}$};
    \node [block, below of=B, yshift=-1.8cm,minimum width=3.0cm, minimum height=1cm] (Blam) {Problem $P_{\lambda}$};
    \draw [-latex] (C.south) -- node[pos=0.5,right,yshift=0cm] {\scriptsize\parbox{2cm}{KKT}} (Clam.north);
    \draw [-latex] (B.south) -- node[pos=0.5,right] {\scriptsize\parbox{1.5cm}{Minimum Principle}} (Blam.north);
    \draw [-latex] (Clam.west) -- node[pos=0.5,below] {\scriptsize{$N \to \infty$}} node[pos=0.5,above] {\scriptsize{\bf Covector mapping, Eq. \eqref{eq:mapping}}} (Blam.east);
    \node (Previous) [dashed,rounded corners,draw=black,fit={(B) (C)}, minimum height = 1.6cm] {};
    \node[above left, inner sep=2pt] at (Previous.north east) {\scriptsize{\bf Previous work, see \cite{cichella2020optimal}}};
    \node (Current) [dashed,rounded corners,draw=black,fit={(Clam) (Blam)}, minimum height = 1.6cm] {};
    \node[below left, inner sep=2pt] at (Current.south east) {\scriptsize{\bf This work}};
\end{tikzpicture}%
%}
\caption{Diagram of the covector mapping principle for Bernstein approximation. The solution to Problem  $P_{N \lambda}$ converges to that of Problem $P_\lambda$ as $N \to \infty$.}
\label{fig:covectormapping}
\end{figure}
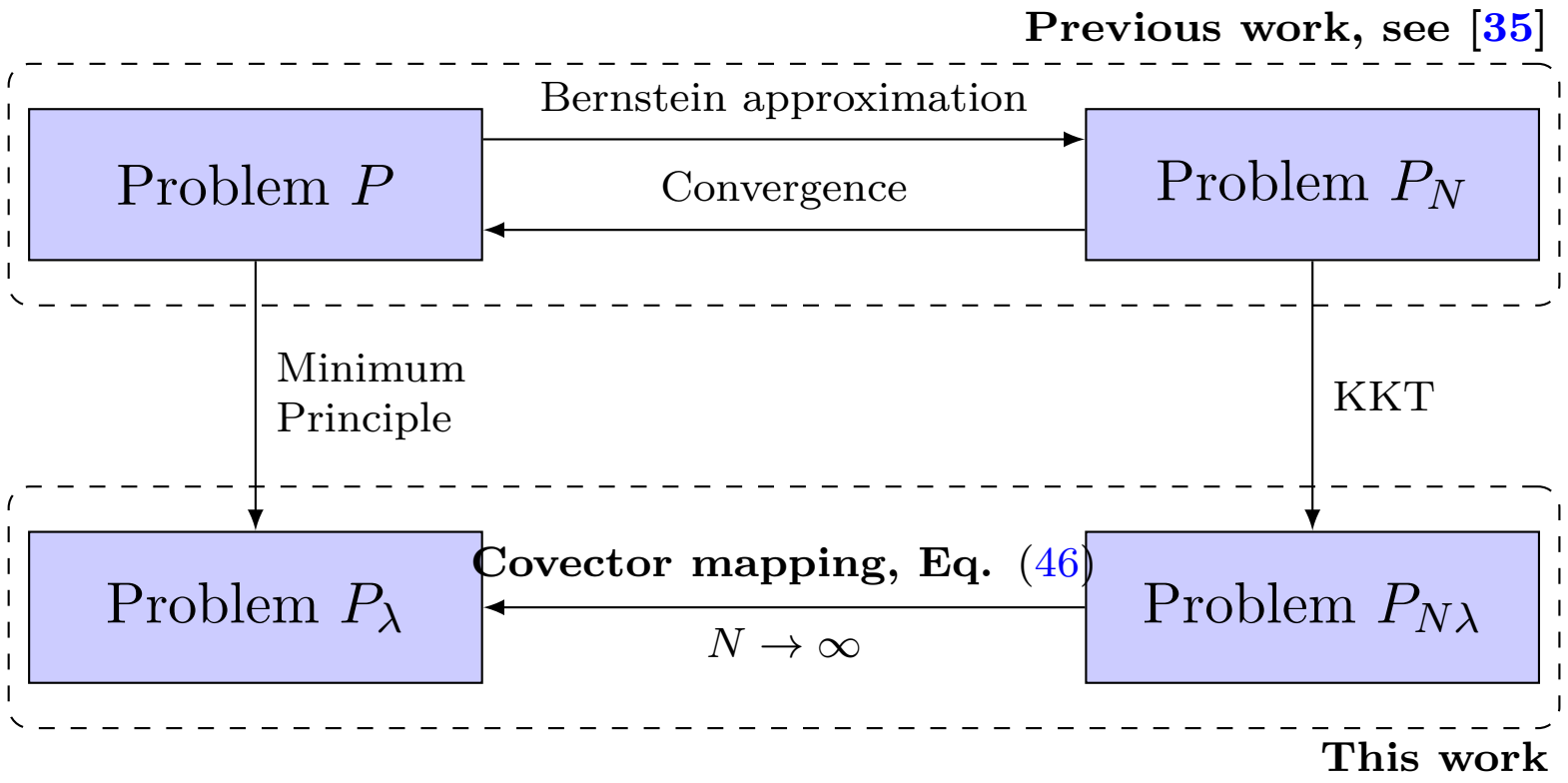

\section{Numerical examples} \label{sec:numericalresults}

\subsection{Example 1: nonlinear one-dimensional problem}
Consider the following problem, taken from \cite{garg2011direct}: determine $y:[0,5]\to\IR$ and $u:[0,5]\to\IR$ that minimize
\begin{equation} \label{eq:ex1cost}
J = \frac{1}{2}\int_0^5 \left( y(t) + u^2(t) \right) dt
\end{equation}
subject to
\begin{equation} \label{eq:ex1dyn}
\dot{y}(t) = 2y(t) + 2u(t)\sqrt{y(t)} \, , \quad \forall t\in[0,5] \, , \qquad y(0) = 2 \, , \quad y(5) = 1 \, .
\end{equation}
The exact solution is given by \cite{garg2011direct}
\begin{equation} \label{eq:ex1exact}
y^*(t) = z^2(t) \, , \qquad u^*(t) = -\lambda_z(t) \, , \qquad \lambda_y^*(t) = \frac{\lambda_z(t)}{2z(t)} \, ,
\end{equation}
where
$$
\begin{bmatrix} z(t) \\ \lambda_z(t) \end{bmatrix} = \exp(\bm{A}t)\begin{bmatrix} \sqrt{2} \\ \lambda_{z0} \end{bmatrix} \, , \qquad \bm{A} = \begin{bmatrix} 1 & -1 \\ -1 & -1 \end{bmatrix} \, , \qquad \lambda_{z0} = \frac{1-\sqrt{2}\,B_{11}}{B_{12}} \, ,
$$
and $B_{11}$, $B_{12}$ are the entries in the first row of $\exp(5\bm{A})$. The expression \eqref{eq:ex1exact} follows from the change of variables $z=\sqrt{y}$, which transforms the problem into a linear quadratic one with $\dot{z} = z + u$ and running cost $\frac{1}{2}(z^2+u^2)$.

To cast the problem in the form of Problem $P$, time is normalized as $t = 5\tau$, $\tau\in[0,1]$, so that the dynamics and the running cost are multiplied by $5$. With this normalization, the costate of the normalized problem at $\tau$ coincides with $\lambda_y^*(5\tau)$. Problem $P_N$ is solved with $\delta_P^N = 0$, i.e., with the dynamic constraint \eqref{eq:discdynamicconstraint} imposed as an equality at the Bernstein coefficients, and with the boundary conditions imposed on the first and last coefficients. The NLP is formulated in CasADi \cite{andersson2019casadi} and solved with IPOPT \cite{wachter2006implementation} with tolerance $10^{-12}$, for $N$ between $15$ and $320$. Since $\sqrt{y}$ is not differentiable at $y=0$ and the exact state reaches a minimum of approximately $4.8\cdot10^{-3}$ near $t\approx 2.6$, the safeguard $\bar{y}_{j,N}\geq 10^{-6}$ is added to the NLP; it is inactive at all the solutions reported below, so that these are solutions of Problem $P_N$. For $N\leq 10$, the solver returned points at which the safeguard is active, and these values of $N$ are not reported. For every $N$, the solver is initialized with $\bar{y}_{j,N} = 2 - \frac{j}{N}$ (the linear interpolation of the boundary conditions) and $\bar{u}_{j,N} = 0$, i.e., without using the exact solution. The costate is estimated from the KKT multipliers returned by the solver using the covector mapping \eqref{eq:mapping}, i.e., $\lambda_y(t_j)\approx(N+1)\bar{\lambda}_{j,N}$, with $t_j = 5j/N$. No closure conditions are imposed.

\begin{figure}[h]
	\centering
	\includegraphics[width = 1\columnwidth]{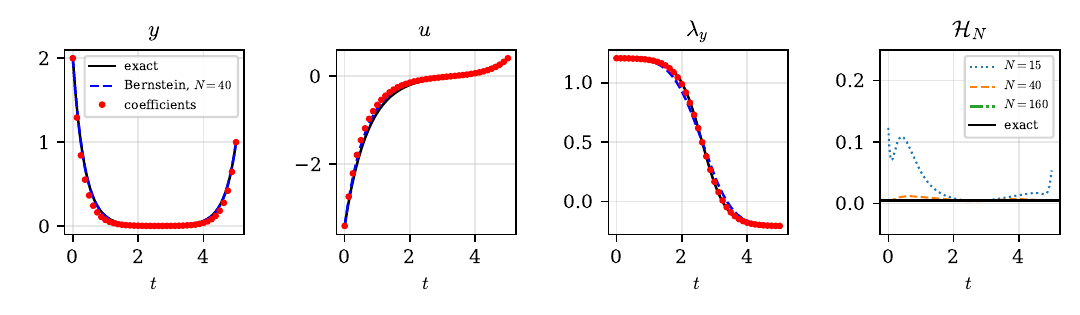}
	\caption{Example 1, $N=40$: state, control and costate (Bernstein polynomials, dashed; Bernstein coefficients, dots; exact solution, solid), and Hamiltonian approximation $\mathcal{H}_N(t)$ for $N=15,40,160$.}
	\label{fig:ex1solution}
\end{figure}

Figure \ref{fig:ex1solution} shows the solution for $N=40$ together with the exact solution. The Bernstein polynomials and their coefficients follow the exact state, control and costate, including the steep boundary layers of the state and of the control near $t=0$ and $t=5$. The right panel of Figure \ref{fig:ex1solution} shows the Hamiltonian approximation $\mathcal{H}_N(t)$ of Remark \ref{rem:hamiltonian}, computed from the Bernstein polynomials, for $N=15$, $40$ and $160$. Since the problem is autonomous, the exact Hamiltonian is constant ($\mathcal{H}\approx 4.8\cdot 10^{-3}$), and $\mathcal{H}_N(t)$ approaches this constant as $N$ increases.

\begin{figure}[h]
	\centering
	\includegraphics[width = 1\columnwidth]{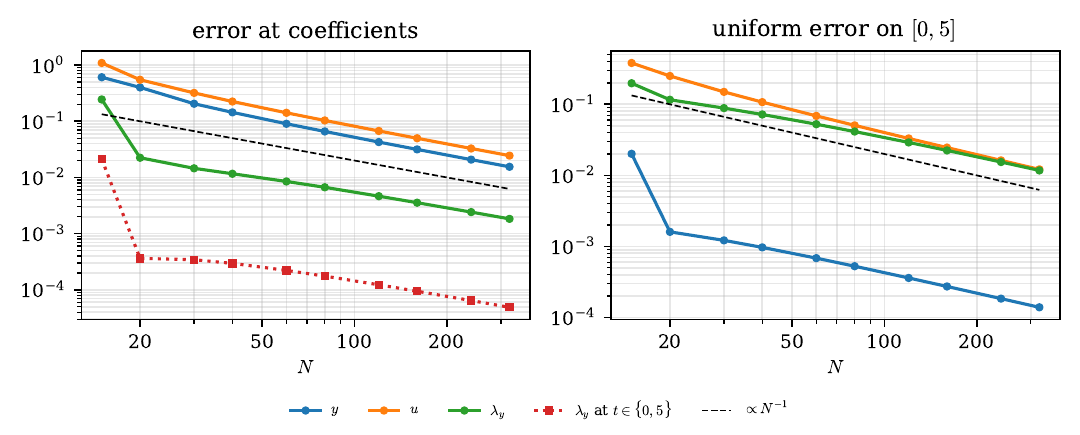}
	\caption{Example 1: maximum errors at the Bernstein coefficients (left) and maximum errors of the Bernstein polynomials over $[0,5]$ (right), as functions of $N$. The dashed line is proportional to $N^{-1}$.}
	\label{fig:ex1errors}
\end{figure}

Figure \ref{fig:ex1errors} reports, as functions of $N$, the maximum errors at the Bernstein coefficients, $\max_j|\bar{y}_{j,N}-y^*(t_j)|$, $\max_j|\bar{u}_{j,N}-u^*(t_j)|$ and $\max_j|(N+1)\bar{\lambda}_{j,N}-\lambda_y^*(t_j)|$ (left), and the maximum errors of the Bernstein polynomials over $[0,5]$ (right). After a preasymptotic range ($N< 20$), all errors decrease at the rate $N^{-1}$, which is consistent with Corollary \ref{cor:c2_dual} and with the saturation of Bernstein approximation at order $N^{-1}$. For $N=320$, the costate error at the coefficients is $1.8\cdot10^{-3}$. The left panel also shows the costate error at the end points, $\max\{|(N+1)\bar{\lambda}_{0,N}-\lambda_y^*(0)|,|(N+1)\bar{\lambda}_{N,N}-\lambda_y^*(5)|\}$, which equals $4.9\cdot10^{-5}$ for $N=320$. This confirms numerically that the transversality conditions are recovered from the KKT conditions without closure conditions (see Equations \eqref{eq:scaledzero} and \eqref{eq:scaledN}).


The multiplier $\bar{\bm\nu}$ of the boundary conditions provides the sensitivity of the optimal cost to the boundary data. Writing the boundary conditions as $e_1 = y(0)-a=0$ and $e_2 = y(5)-b=0$, with $a=2$ and $b=1$, Equations \eqref{eq:Plambda3} and \eqref{eq:Plambda4} give $\nu_1 = -\lambda_y^*(0)$ and $\nu_2 = \lambda_y^*(5)$, and the optimal cost $J^*(a,b)$ satisfies $\partial J^*/\partial a = \lambda_y^*(0) = -\nu_1$ and $\partial J^*/\partial b = -\lambda_y^*(5) = -\nu_2$. Table \ref{tab:ex1sens} compares $-\bar{\nu}_1$ and $-\bar{\nu}_2$, obtained from the solution of Problem $P_N$, with central finite differences (step $10^{-4}$) of the exact optimal cost $J^*(a,b)$, which is computed by quadrature of the exact solution and does not use the costate. The estimates agree with the finite differences to within $1.2\cdot10^{-4}$ for $N=320$, with errors decreasing as $N^{-1}$. Hence, a single solution of Problem $P_N$ predicts, to first order, the change of the optimal cost under a perturbation of the initial or final state, without solving additional problems.

\begin{table}[h]
\centering
\caption{Example 1: sensitivity of the optimal cost to the boundary values $a=y(0)$ and $b=y(5)$. Exact values are central finite differences of $J^*(a,b)$: $\partial J^*/\partial a = 1.206258$, $\partial J^*/\partial b = 0.205409$.}
\label{tab:ex1sens}
\begin{tabular}{rcccc}
\hline
$N$ & $-\bar{\nu}_1$ & $\vert -\bar{\nu}_1 - \partial J^*/\partial a\vert$ & $-\bar{\nu}_2$ & $\vert -\bar{\nu}_2 - \partial J^*/\partial b\vert$ \\
\hline
20 & 1.206835 & $5.8\cdot10^{-4}$ & 0.206579 & $1.2\cdot10^{-3}$ \\
40 & 1.206631 & $3.7\cdot10^{-4}$ & 0.206171 & $7.6\cdot10^{-4}$ \\
80 & 1.206471 & $2.1\cdot10^{-4}$ & 0.205846 & $4.4\cdot10^{-4}$ \\
160 & 1.206372 & $1.1\cdot10^{-4}$ & 0.205643 & $2.3\cdot10^{-4}$ \\
320 & 1.206317 & $5.9\cdot10^{-5}$ & 0.205530 & $1.2\cdot10^{-4}$ \\
\hline
\end{tabular}
\end{table}

\subsection{Example 2: problem with an active mixed state-input constraint}
The second example tests the estimation of the multiplier $\bm{\mu}$ of the inequality constraints, which is not exercised by Example 1. Consider the problem: determine $x:[0,1]\to\IR$ and $u:[0,1]\to\IR$ that minimize
\begin{equation} \label{eq:ex2cost}
J = \frac{1}{2}\int_0^1 \left( x^2(t) + u^2(t) \right) dt
\end{equation}
subject to
\begin{equation} \label{eq:ex2dyn}
\dot{x}(t) = u(t) \, , \qquad x(0) = 1 \, , \qquad h(x(t),u(t)) = -u(t) - c\,x(t) \leq 0 \, , \quad \forall t\in[0,1] \, ,
\end{equation}
with $c=\frac{1}{2}$. The constraint $u \geq -c\,x$ bounds the rate of decrease of the state relative to the state itself, and it depends on both $x$ and $u$; the final state is free. Without the constraint, the optimal control is $u(t) = -x(0)\sinh(1-t)/\cosh(1)$, which satisfies $u/x = -\tanh(1-t) < -c$ near $t=0$, so the constraint is active on an initial arc. Since the cost is convex, the dynamics are linear and the constraint is linear, the first order conditions in Problem $P_\lambda$ are also sufficient, and the solution is obtained by solving them on each arc and joining the arcs continuously. With $t_1 = 1-\operatorname{artanh}(c) \approx 0.4507$, the exact solution is
\begin{equation} \label{eq:ex2exact}
\begin{split}
& t\in[0,t_1]: \quad x^* = e^{-ct} \, , \quad u^* = -c\,x^* \, , \quad \lambda^* = k\,e^{-ct} + (c-k)\,e^{c(t-2t_1)} \, , \quad \mu^* = \lambda^* - c\,x^* \, , \\
& t\in[t_1,1]: \quad x^* = \alpha\cosh(1-t) \, , \quad u^* = -\alpha\sinh(1-t) \, , \quad \lambda^* = -u^* \, , \quad \mu^* = 0 \, ,
\end{split}
\end{equation}
where $k = \frac{1+c^2}{2c}$ and $\alpha = e^{-ct_1}/\cosh(1-t_1)$. On $[0,t_1]$, $\lambda^*$ solves $\dot{\lambda} = -x + c\,\mu$ with $\mu = u + \lambda$ from $\mathcal{L}_u=0$, and $\mu^*\geq0$, with $\mu^*(t_1)=0$; on $[t_1,1]$, the transversality condition $\lambda^*(1)=0$ holds. The control $u^*$ and the multiplier $\mu^*$ are continuous with a corner at $t_1$, and the state and the costate are of class $\mathcal{C}^1$, with a corner of $\dot{x}^*$ and $\dot{\lambda}^*$ at $t_1$. Hence, Assumption \ref{asm:plambdaexists} holds, whereas the regularity required by Corollary \ref{cor:c2_dual} does not. Since the problem is invariant under the scaling $(x,u)\mapsto(sx,su)$, $s>0$, the optimal cost is homogeneous of degree two in $x(0)>0$, and $\partial J^*/\partial x(0) = 2J^*/x(0) = \lambda^*(0) \approx 0.772111$.

Problem $P_N$ is solved as in Example 1, with $\delta_P^N=0$, the constraint $h(\bar{x}_{j,N},\bar{u}_{j,N})\leq0$ imposed at each pair of Bernstein coefficients, the initial guess $\bar{x}_{j,N}=1$, $\bar{u}_{j,N}=0$, and without relaxation of the inequality constraints by the solver. The costate and the constraint multiplier are estimated by the covector mapping \eqref{eq:mapping}, $\lambda(t_j)\approx(N+1)\bar{\lambda}_{j,N}$ and $\mu(t_j)\approx(N+1)\bar{\mu}_{j,N}$, with $t_j = j/N$.

\begin{figure}[h]
	\centering
	\includegraphics[width = 1\columnwidth]{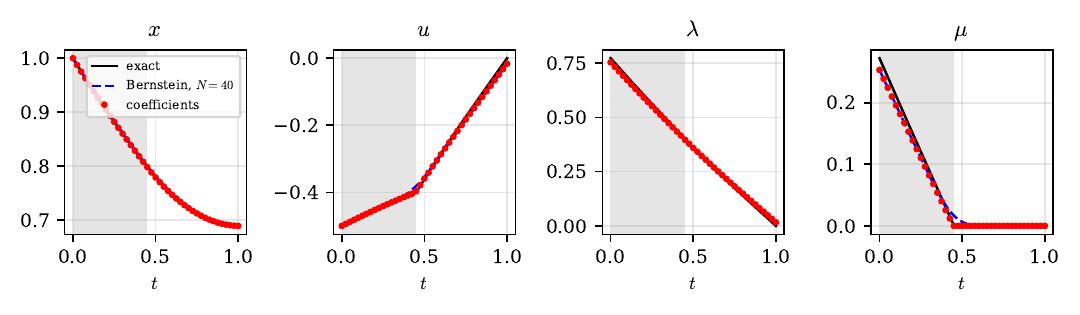}
	\caption{Example 2, $N=40$: state, control, costate and constraint multiplier (Bernstein polynomials, dashed; Bernstein coefficients, dots; exact solution, solid). The shaded region is the constrained arc $[0,t_1]$.}
	\label{fig:ex2solution}
\end{figure}

\begin{figure}[h]
	\centering
	\includegraphics[width = 1\columnwidth]{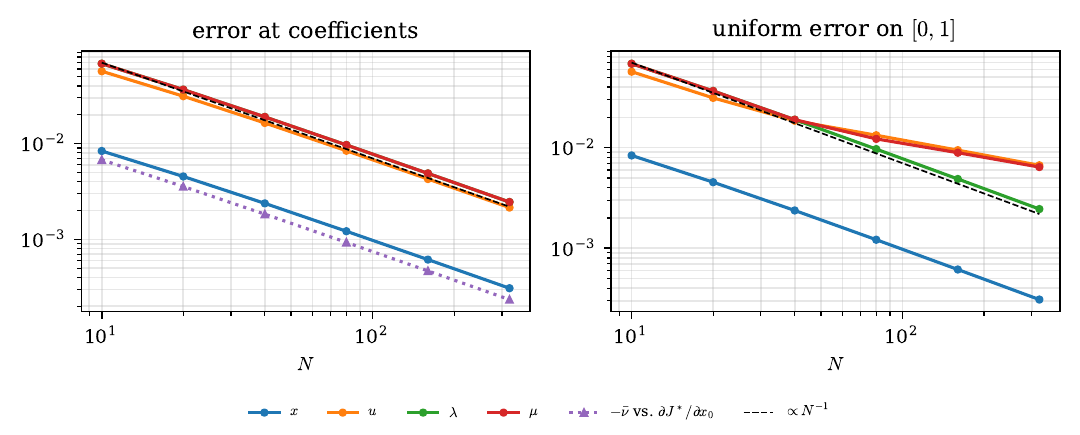}
	\caption{Example 2: maximum errors at the Bernstein coefficients (left) and maximum errors of the Bernstein polynomials over $[0,1]$ (right), as functions of $N$, and error of $-\bar{\nu}$ as an estimate of $\partial J^*/\partial x(0)$. The dashed line is proportional to $N^{-1}$.}
	\label{fig:ex2errors}
\end{figure}

Figure \ref{fig:ex2solution} shows the solution for $N=40$. The estimated multiplier $(N+1)\bar{\mu}_{j,N}$ is positive on the active arc and below $3\cdot10^{-8}$ on the inactive arc, so that the active set of Problem $P_N$ identifies the constrained arc; the junction time estimated as the midpoint between the last coefficient at which the constraint is active, i.e., $h(\bar{x}_{j,N},\bar{u}_{j,N})>-10^{-6}$, and the next one is within $1/N$ of $t_1$ for all $N$ considered (e.g., $8.7\cdot10^{-4}$ for $N=320$). Figure \ref{fig:ex2errors} reports the errors as functions of $N$. At the Bernstein coefficients, the errors of the state, the control, the costate and the multiplier $\mu$ decrease as $N^{-1}$, a rate that is observed here but not guaranteed by Corollary \ref{cor:c2_dual}; for $N=320$ they are $3.1\cdot10^{-4}$, $2.1\cdot10^{-3}$, $2.5\cdot10^{-3}$ and $2.5\cdot10^{-3}$, respectively. The uniform errors of $x_N(t)$ and $\lambda_N(t)/w$ also decrease as $N^{-1}$, whereas the uniform errors of $u_N(t)$ and $\mu_N(t)/w$ decrease more slowly, approximately as $N^{-1/2}$, because of the corner at $t_1$; this is consistent with the bound $W(N^{-1/2})$ of Lemma \ref{lem:bernsteinapprox}, which, for a Lipschitz function, is of order $N^{-1/2}$. The error of $-\bar{\nu}$ as an estimate of $\partial J^*/\partial x(0)$ is $2.3\cdot10^{-4}$ for $N=320$. Finally, since $h$ is linear, by Remark \ref{rem:relaxation} the constraint is satisfied by the Bernstein polynomials for all $t\in[0,1]$, and not only at the coefficients: for all $N$ considered, $\max_{t\in[0,1]} h(x_N(t),u_N(t)) \leq 0$ on a grid of $4001$ points.

\section{Conclusions} \label{sec:conclusions}
This paper proposed a numerical method for costate estimation of nonlinear constrained optimal control problems using Bernstein polynomials, in which the dynamic and path constraints are imposed at the Bernstein coefficients.
A rigorous analysis is provided that shows convergence of the costate estimates to the dual variables of the continuous-time problem. To this end, a set of conditions are derived under which the Karush Kuhn Tucker multipliers of the NLP converge to the costates of the optimal control problem. This led to the formulation of the Covector Mapping Theorem for Bernstein approximation. Owing to the structure of the Bernstein differentiation matrix, the discrete transversality conditions are contained in the KKT conditions, and no additional closure conditions need to be imposed.
  The theoretical findings are illustrated by two numerical examples with known exact solutions, one of which has an active mixed state-input constraint, and the multipliers of the boundary conditions are shown to predict the sensitivity of the optimal cost to the boundary data.

\bibliographystyle{unsrt}
\bibliography{refs}

\end{document}